\documentclass[reqno,oneside,11pt,draft]{amsart}
\usepackage{amssymb,epsfig,verbatim}
\addtocounter{section}{-1}             
\numberwithin{equation}{section}       
\theoremstyle{plain}
\newtheorem{Thm}{Theorem}[section]
\newtheorem{Prop}[Thm]{Proposition}
\newtheorem{Lemma}[Thm]{Lemma}
\newtheorem{Cor}[Thm]{Corollary}
\newtheorem*{ThmA}{Theorem A}
\newtheorem*{ThmA'}{Theorem A'}
\newtheorem*{CorB}{Corollary B}
\newtheorem*{CorC}{Corollary C}
\newtheorem*{ConjD}{Conjecture D}

\theoremstyle{definition}
\newtheorem{Example}[Thm]{Example}
\newtheorem{Def}[Thm]{Definition}
\newtheorem{Remark}[Thm]{Remark}
\newtheorem*{Ackn}{Acknowledgement}

\newcommand{\C}{{\mathbf{C}}}
\newcommand{\D}{{\mathbf{D}}}
\newcommand{\N}{{\mathbf{N}}}

\newcommand{\R}{{\mathbf{R}}}

\newcommand{\CC}[1]{{\mathbf{C}^{#1}}}
\newcommand{\PP}[1]{{\mathbf{P}^{#1}}}

\newcommand{\cC}{{\mathcal{C}}}
\newcommand{\cE}{{\mathcal{E}}}

\newcommand{\cH}{{\mathcal{H}}}
\newcommand{\cO}{{\mathcal{O}}}
\newcommand{\cS}{{\mathcal{S}}}

\newcommand{\wmu}{\widehat{\mu}}
\newcommand{\wphi}{\widetilde{\phi}}
\renewcommand{\=}{:=}
\renewcommand{\a}{\alpha}
\newcommand{\e}{\varepsilon}
\newcommand{\vp}{\varphi}
\renewcommand{\l}{\lambda}

\newcommand{\area}{\mathrm{Area\, }}
\newcommand{\dist}{\mathrm{dist}}

\newcommand{\sing}{\mathrm{sing}\,}
\newcommand{\vol}{\mathrm{Vol\, }}

\newcommand{\self}{\circlearrowleft}

\newcommand{\notto}{\nrightarrow}
\newcommand{\Lloc}{L^1_{\mathrm{loc}}}
\newcommand{\FS}{Forn{\ae}ss and Sibony }

\begin{document}
%
%
\title{Brolin's Theorem for curves in two complex dimensions}

\date{\today}
\author{Charles Favre \and Mattias Jonsson}
\address{Universit\'e Paris 7\\
         UFR de Math\'ematiques\\
         Equipe G\'eom\'etrie et Dynamique\\
         F-75251 Paris Cedex 05\\
         France}
\email{favre@math.jussieu.fr}
\address{Department of Mathematics\\ 
         University of Michigan\\
         Ann Arbor, MI 48109-1109\\
         USA}
\email{mattiasj@umich.edu}
\thanks{Second author supported by the Swedish Foundation 
  for Cooperation in Research and Higher Education (STINT)}
%
%
\begin{abstract}
  Given a holomorphic mapping $f:\PP{2}\self$ of degree $d\ge2$
  we give sufficient conditions on a positive closed $(1,1)$
  current $S$ of unit mass under which 
  $d^{-n}f^{n*}S$ converges to the Green current as $n\to\infty$.
\end{abstract}

%
%
\maketitle
%
%

\section{Introduction}\label{intro}

In 1965 H.~Brolin~\cite{Bro1} proved a remarkable result about 
the distribution 
of preimages of points for polynomial maps in one variable: if
$f(z)=z^d+\dots$ is a polynomial of degree $d\ge2$, then there is
a set $\cE$ with $\#\cE\le1$ such that if $a\in\C\setminus\cE$, then
\begin{equation}\label{Me1}
  \frac1{d^n}\sum_{f^nz=a}\delta_z\to\mu\quad\text{as $n\to\infty$},
\end{equation}
where $\mu$ is harmonic measure on the filled Julia set of $f$. 
In particular, the limit in~\eqref{Me1} is independent of $a$. Further,
the \emph{exceptional set} $\cE$ is empty unless $f$ is affinely
conjugate to $z\mapsto z^d$, in which case $\cE=\{0\}$, the 
totally invariant point.

Lyubich~\cite{Lyu1} and Freire, Lopez and Ma{\~n}{\'e}~\cite{FrLoMa1}
later generalized Brolin's theorem to rational maps of the 
Riemann sphere $\PP{1}$, with $\#\cE\le2$.

In this paper we prove a version of Brolin's Theorem in two complex
dimensions.

\begin{ThmA}
  Let $f:\PP{2}\self$ be a holomorphic mapping of algebraic 
  degree $d\ge2$. Then there is a totally invariant,
  algebraic set $\cE_1$ consisting of at most three complex lines
  and a finite, totally invariant set $\cE_2$ with the following
  property:
  If $S$ is a positive closed current on $\PP{2}$ of bidegree
  $(1,1)$ and mass $1$ such that
  \begin{itemize}
  \item[(i)]
    $S$ does not charge any irreducible
    component of $\cE_1$;
  \item[(ii)]
    $S$ has a bounded local potential at each point of $\cE_2$;
  \end{itemize}  
  then we have the convergence 
  \begin{equation}\label{Me2}
    \frac1{d^n}f^{n*}S\to T\quad\text{as $n\to\infty$}.
  \end{equation}
\end{ThmA} 

Here $T$ is the Green current of $f$, defined as the limit
(in the sense of currents) of $d^{-n}f^{n*}\omega$ as $n\to\infty$, 
where $\omega$ is the Fubini-Study form on $\PP{2}$. See 
Section~\ref{background} for more details.

As a consequence we have the following result on the distribution of
the preimages of curves. The space of curves in $\PP{2}$ of
degree $k$ may be identified with $\PP{N}$ for some $N=N(k)$.
\begin{CorB}
  Let $f$ be as in Theorem~A and let $k\ge1$. 
  Let $\cE^*$ be the set of curves $C\in\PP{N}$ such that
  \begin{equation*}
    \frac1{d^nk}f^{n*}[C]\notto T\quad\text{as $n\to\infty$}.
  \end{equation*}
  Then $\cE^*$ is an algebraic proper subset of $\PP{N}$.
\end{CorB}

In a similar way, the space $\mathrm{Hol}_d$ of holomorphic maps of $\PP{2}$
of degree $d$ can be identified with a Zariski open set of some
$\PP{M}$.
\begin{CorC}
  There exists an algebraic proper subset ${\mathcal H} \subset
  \mathrm{Hol}_d$ such that for any $f \notin {\mathcal H}$ the convergence
  \begin{equation*}
    \frac1{d^n}f^{n*}S\to T\quad\text{as $n\to\infty$}
  \end{equation*}
  holds for all positive closed $(1,1)$ currents $S$ of mass $1$.
\end{CorC}

The set $\cE_2$ is a little mysterious but contains the following
two types of points:
\begin{itemize}
\item[($\alpha$)]
  totally invariant points on totally invariant curves;
\item[($\beta$)]
  homogeneous points, that is $f$ preserves the pencil of lines passing
  through the point.
\end{itemize}
It can be shown the set of such points contains at most three
elements.  We postpone to a later study the fact that these are the
only types of points in $\cE_2$.

Denote by $\nu(p,S)$ the Lelong number of the positive closed
current $S$ at $p$. We then have
\begin{ThmA'}
  Let $f, S$ be as in Theorem~A, and assume $\cE_2$ is
  reduced to points of type $(\alpha)$ and $(\beta)$. Then the following
  statements are equivalent:
  \begin{itemize}
  \item[(1)]
     $d^{-n}f^{n*} S\to T$ as $n\to\infty$;
  \item[(2)]$S$ does not charge any irreducible
    component of $\cE_1$ and 
    $\nu(p,S)=0$ at each point $p\in\cE_2$.
  \end{itemize}  
\end{ThmA'} 

We mention the following
\begin{ConjD}
 Let $f, S$ be as in Theorem~A. Then $d^{-n}f^{n*} S$ converges to a
 current $\overline{S}$ satisfying the invariance relation
 $f^*\overline{S}= d \overline{S}$.
\end{ConjD}

Results in the direction of our paper were previously obtained by
\FS\cite{FoSi2}, who proved a weaker form of Corollary~C
with $\mathcal H$ a \emph{countable} union of algebraic sets in
$\mathrm{Hol}_d$.
Later Russakovskii and Shiffman~\cite{RuSh1} proved a version of
Corollary~B: for any holomorphic mapping $f:\PP{2}\self$ there exists
a pluripolar set $\cE^*\subset\PP{2*}$ such that if $L$ is a line
in $\PP{2*}\setminus\cE^*$, then $d^{-n}f^{n*}[L]\to T$ as
$n\to\infty$. In fact, their result applies also to certain
rational maps of $\PP{k}$, $k\ge 2$ and pullbacks of planes of higher 
codimension.

Other related results include Briend and Duval~\cite{BrDu2}, who 
recently proved Brolin's theorem for
preimages of \emph{points} under holomorphic maps of $\mathbf{P}^k$.
A version of Theorem~A was proven by 
\FS~\cite{FoSi7} for H\'enon maps (see also \cite{BeSm9}) 
and by Favre and Guedj~\cite{FaGu1} for birational maps
(see also~\cite{Fav4}). 

The main ingredient in our proofs of Theorems~A and~A', 
as well as in most previous approaches, are 
\emph{volume estimates}. They come in two forms.

The first type of volume estimates are dynamical and aim at bounding
$\vol f^nE$ from below for any Borel set $E$. Such estimates are
related to the rate of recurrence of the critical set. 
In previous work, assumptions on the dynamics were made to get the
required volume estimates. A main novelty of this paper is that we 
are able to control volume decay for arbitrary holomorphic maps.

More precisely we show that the phase space $\PP{2}$ splits
naturally into two parts: the exceptional set $\cE=\cE_1\cup\cE_2$ 
and its complement.
Outside $\cE$, the critical set is not too recurrent and 
$\vol f^nE$ does not decay too fast. Near $\cE$, on the other hand,
$\vol f^nE$ may a priori decay quite rapidly, but precise information on
the structure of $\cE$ allows a good understanding of the dynamics
and in particular of volume decay.

To obtain this partition we first study asymptotic
volume decay along orbits and relate it to the growth of two
algebraic quantities: the multiplicity of the vanishing of the
Jacobian determinant, and the generic rate of contraction. A key
result in the paper is the understanding of the asymptotic behavior of
these multiplicities under iteration. In particular, we characterize
the locus where these asymptotic multiplicities are maximal,
giving rise to the exceptional sets $\cE_1$ and $\cE_2$ in Theorem~A.
Semicontinuity properties of the multiplicities imply that these sets
are algebraic and present strong recurrence properties: they are in 
fact totally invariant.

The second type of volume estimates involve pluripotential theory. 
We estimate the volume of sublevel sets of plurisubharmonic (psh)
functions using the Kiselman-Skoda theorem:
asymptotics of these volumes are small exactly when Lelong numbers are
small. 

In this context, to show that certain Lelong numbers decay under
iteration, we make use of Kiselman numbers (or directional Lelong
numbers).  These allow us to deduce dynamical information in a
neighborhood of an invariant curve from the dynamics on the curve
itself. We believe this technique could prove useful in other
situations, too.

We also believe our result to be true in any dimension but the
description of $\cE$ and hence the control of decay of volumes around
$\cE$ become much harder than in dimension $2$.

The organization of this paper is as follows.  We briefly recall some
facts from holomorphic dynamics and reduce Theorem~A to an estimate of
the size of images of balls in Section~\ref{background}.  In
Section~\ref{pluripot}, we state some pluripotential facts that we use
in the paper. In particular we investigate the behavior of Kiselman
numbers as one weight degenerates. Same results appeared independently
in~\cite{mim}. In Section~\ref{multi} we define three asymptotic
multiplicities related to volume decay.  These multiplicities are used
to define the exceptional sets $\cE_1$ and $\cE_2$, and we study the
latter sets in Sections~\ref{exc1} and~\ref{exc2}.  The next two
sections are devoted to volume estimates outside $\cE_1\cup\cE_2$
(Section~\ref{volume0}) and near $\cE_1\setminus\cE_2$
(Section~\ref{volume1}).  In Section~\ref{lelong}, we show a useful
technical result about Lelong numbers of pull-backs of currents near a
totally invariant curve.  After these estimates, we prove Theorem~A
and Corollaries~B and~C in Section~\ref{proofmain}. Finally, we prove
Theorem~A' and discuss the existence of totally invariant currents in
Section~\ref{misc}.

\begin{Ackn}
  This paper was partially written when the authors were visiting 
  IMPA, Rio de Janeiro, and they wish to thank the department
  for its hospitality and support. 
\end{Ackn}

%
%
%
%
\section{Background and Reduction}\label{background}

In this section we recall some known facts about holomorphic 
mappings of $\PP{2}$; see e.g.~\cite{FoSi4} for more information.
We also reduce the proof of Theorem~A to an estimate
on the sizes of images of balls.

Let $f:\PP{2}\self$ be a holomorphic map of degree $d\ge2$. 
This means that $f=[P(z,w,t):Q(z,w,t):R(z,w,t)]$, 
where $P$, $Q$, and $R$ are homogeneous polynomials of degree
$d$ with no nontrivial common zero. 

Let $T$ be a positive closed $(1,1)$ current on $\PP{2}$ and take a
local plurisubharmonic (psh) 
potential $T= dd^c u$ around $p \in \PP{2}$. One defines
locally at any point in $f^{-1}\{p\}$ the positive closed $(1,1)$ current
$f^* T \= dd^c (u \circ f)$. This does not depend on the choice of $u$
and induces a continuous linear operator on the set of positive
closed $(1,1)$ currents on $\PP{2}$.
One can project $f^*$  to an action $f^*$ on $H^{1,1}_{\R}(\PP{2})
\simeq \R$. This latter is given by the multiplication by $d$.

Let $\omega$ be the  Fubini-Study K\"ahler form on $\PP{2}$. The positive
closed currents $f^*\omega$ and $d\,\omega$ are cohomologous, 
one can hence find a continuous function $u$ 
such that $f^*\omega=d\,\omega+dd^c u$.
Iterating this equality $n$ times yields 
$d^{-n}f^{n*}\omega=\omega+dd^c(\sum_{j=1}^nd^{-j}u\circ f^{j-1})$. 
This latter series converges uniformly on $\PP{2}$ to a continuous 
function $G$ and one finally infers 
\begin{equation}\label{Me3}
  \frac1{d^n}f^{n*}\omega\to T\quad\text{as $n\to\infty$},
\end{equation}
where $T\=\omega+dd^cG$ is called the \emph{Green current} of
$f$. It satisfies the invariance property $d^{-1}f^*T=T$.
Replacing $\omega$ in~\eqref{Me3} by a general positive closed current
$S$ of unit mass leads to~\eqref{Me2}; the purpose of this paper is to 
investigate exactly for what currents $S$ this convergence holds.

As stated in Theorem~A, 
the exceptional currents (for which~\eqref{Me2} may fail)
will be connected with totally invariant algebraic sets, and we recall
the following two results.
\begin{Prop}~\label{MP1}~\cite{FoSi1,ShShUe1,Da2,CeLN1}
  Let $f:\PP{2}\self$ be holomorphic of degree $d\ge2$. 
  Then the following holds:
  \begin{itemize}
  \item[(i)]
    any (possibly reducible) totally invariant curve $V\supset f^{-1}V$
    is a union of at most three lines; if there are
    three lines, then they are in generic position;
    further, the set of intersection points between different
    lines is totally invariant;
  \item[(ii)]
    any finite totally invariant set $X\supset f^{-1}X$
    is contained in the singular locus 
    of the critical value set of $f$;
    the cardinality of the union of all such 
    sets $X$ is bounded by a number depending only on $d$.
  \end{itemize}
\end{Prop}

As a start of the proof of Theorem~A we will reduce~\eqref{Me2}
to a study of sizes of images of balls under iterates of $f$.
Namely, suppose that $S$ is a positive, closed current cohomologous
to $\omega$ for
which~\eqref{Me2} fails. Then we may write
\begin{equation*}
  S=\omega+dd^cu,
\end{equation*}
where $u\le0$ is a quasi-plurisubharmonic (qpsh) function on $\PP{2}$.
It then follows that for all $n\ge0$,
\begin{equation*}
  d^{-n}f^{n*}S= d^{-n}f^{n*}\omega+ d^{-n}dd^c(u\circ f^n),
\end{equation*}
so since~\eqref{Me2} fails and $d^{-n}f^{n*}\omega\to T$, we have 
$d^{-n}dd^c(u\circ f^n)\notto 0$ as $n\to\infty$, 
which is equivalent to 
$d^{-n}u\circ f^n\notto 0$ in $\Lloc$.
By Hartog's Lemma (see~\cite{Ho2})
this implies that there is a ball $B\subset\PP{2}$,
a constant $\alpha>0$ and a sequence $n_j\to\infty$ such that
\begin{equation}\label{Me4}
  f^{n_j}B\subset\{u<-\alpha d^{n_j}\}.
\end{equation}
The rest of the proof consists of showing that~\eqref{Me4} is
not possible if $S$ is a current satisfying the hypotheses of
Theorem~A. This is done by estimating the volume
of $f^{n_j}B$ from below (using dynamics) and the volume
of $\{u<-\alpha d^{n_j}\}$ from above (using pluripotential theory).

%
%
%
%
\section{Some pluripotential theory}\label{pluripot}

In this section we discuss some results from pluripotential
theory. First we need the definition of Lelong numbers.

\begin{Def}
  Let $u$ be a a psh function near the origin in $\CC{2}$.
  We define the Lelong number $\nu(0,u)$ of $u$ at the origin
  to be the supremum of $\nu>0$ such that
  \begin{equation*}
    u(\zeta)\le\nu\log\vert\zeta\vert+O(1)\qquad\text{as $\zeta\to0$}.
  \end{equation*}
\end{Def}

This definition is invariant under local biholomorphisms.
We will need an estimate of the volume of sublevel sets for
a psh function in terms of its Lelong numbers. The following result
is due to Kiselman~\cite{Kis1} and 
relies on previous work of Skoda~\cite{Sk1}.

\begin{Thm}\label{vol-t:2}
  Let $U \subset \C ^2$ be an open set, $K$ a compact subset of $U$, 
  and $u$ a psh function on $U$. 
  For any real number $\alpha<2(\sup_K\nu(z,u))^{-1}$ 
  there exists a constant $C_{\alpha}>0$ such that for any $t\ge0$, 
  the estimate
  \begin{equation}\label{vol-e:10}
    \vol\left(K\cap\{u\le-t\}\right)\le C_{\alpha}\exp(-\alpha t).
  \end{equation}
holds.
\end{Thm}

In Section~\ref{lelong} we will need to work with directional Lelong
numbers or \emph{Kiselman numbers}. We refer to~\cite{Kis2}  
or~\cite{De1} for a detailed exposition.

Let $u$ be a psh function in the unit ball $B \subset \C ^2$ endowed
with coordinates $(z,w)$. Fix a weight $(\a_1, \a_2) \in (\R_+^*)^2$.

The Kiselman number of $u$ at the point $p = (0,0)$ with weight
$(\a_1,\a_2)$ is defined as
\begin{equation*}
  \nu(p,u,(\a_1,\a_2)) 
  \=
  \lim_{r\to 0}\frac{\a_1\a_2}{\log r}
  \sup_{\Delta(r^{1/\a_1})\times\Delta(r^{1/\a_2})}u.
\end{equation*}
For $(\a_1,\a_2)=(1,1)$ we recover the usual Lelong number
\begin{equation*}
  \nu(p,u,(1,1))=\nu(p,u).
\end{equation*}
We have the following homogeneity property
\begin{equation}\label{Ec1}
  \nu(p,u,(\l\a_1,\l \a_2))
  =\l\,\nu(p,u,(\a_1,\a_2)),
\end{equation}
for any $\l>0$.
The inequality
\begin{equation}\label{Ec2}
  \nu(p,u,(\a_1,\a_2))
  \ge\min(\a_1,\a_2)\,\nu(p,u) 
\end{equation}
always holds (see \cite{De1}).

We can now state the following
\begin{Prop}\label{kisel}
  Let $u \in PSH (B)$, $u \le 0$ and assume that the positive closed 
  current $S=dd^c u$ does not charge the curve $\{ z = 0 \}$. Then
  \begin{equation*}
    \lim_{\a\to 0}\sup _{p\in\{z=0\}}\nu(p,S,(\a,1))=0.
  \end{equation*}
\end{Prop}

\begin{Remark}
The same types of considerations were made by Souad Khemeri Mimouni
in~\cite{mim}. In fact, she studies more generally the transformation
of Lelong numbers under any sequence of blowing-ups. In our case, we
blow-up only at the intersection point of the exceptional divisor and
the strict transform of the curve $\{ z = 0 \}$.
\end{Remark}

\begin{proof}[Proof of Proposition~\ref{kisel}.]
  We rely on the approximation process described by Demailly in \cite{De2}.
  Introduce the following Hilbert space  
  \begin{equation*}
    \cH(su)\=\{h\in\cO(B)\ \vert\ 
    |h|^2_{su}\=\int_B|h|^2e^{-2su}<+\infty \}.
  \end{equation*}
  Set
  \begin{equation*}
    u_s\=s^{-1}\sup_{|h|_{su}\le1}\log|h|\in PSH(B).
  \end{equation*}
  One checks that $u_s = (2s)^{-1} \log \sum |\sigma_i|^2$ for any
  orthonormal basis $\{ \sigma_i\}$ of $\cH ( s u )$, and that 
  $u_s$ is the logarithm of
  a real analytic function. The following result connects the singularities
  of $u$ and $u_s$.
  \begin{Lemma}\label{small}
    For any point $p\in B$ and all $(\a_1,\a_2) \in (\R_+^*)^2$, one has
    \begin{enumerate}
    \item
      $\nu(p,u_s,(\a_1,\a_2))\ge\nu(p,u,(\a_1,\a_2))-s^{-1}(\a_1+\a_2)$;
    \item
      $\nu(p,u_s)\le\nu(p,u)$.
    \end{enumerate}
  \end{Lemma}

  We then conclude the proof of the theorem as follows. As $dd^c u$ does
  not charge $\{ z=0 \}$ there exists a point $p = (0,w)$ such that 
  $\nu(p,u)=0$ (by Siu's theorem). In particular by Lemma~\ref{small}~(2), 
  $u_s$ does not charge 
  $\{z=0\}$ either, hence there exists a holomorphic function 
  $h\in\cH(su)$ which does not vanish identically on $\{ z=0\}$. For such a
  function, we apply {\L}ojasiewic's inequality~\cite[p. 243]{Lo1} and get
  \begin{equation*}
    |h(z,w)|+|z|\ge C\dist((z,w),I)^{\theta},
  \end{equation*}
  with $I=h^{-1}(0)\cap\{z=0\}$ and for some constants $C,\theta>0$.
  We infer for any $w$ and any $\a < \theta^{-1}$,
  \begin{equation*}
    \sup_{\Delta( r ^{1/\a})\times w+\Delta(r)}\log|h|
    \ge C'\theta\log r.
  \end{equation*}
  From this, it is easy to see that 
  $\lim_{\a \to 0} \sup_{p\in\{z=0\}}\nu(p,\log|h|,(\a,1))=0$; hence
  \begin{equation*}
   \lim_{\a \to 0}\sup_{p\in\{z=0 \}}\nu(p,u_s,(\a,1))
   \le
   \lim_{\a\to 0}\sup_{p\in\{z=0\}}s^{-1}\nu(p,\log|h|,(\a,1))=0.
 \end{equation*}
 In particular, we get for any $s\ge0$, 
 \begin{equation*}
   \lim_{\a\to 0}\sup_{p\in\{z=0\}}\nu(p,u,(\a,1))
   \le s^{-1}+\lim_{\a\to 0}\sup_{p\in\{z=0\}}\nu(p,u_s,(\a,1))
   \le s^{-1},
 \end{equation*}
 which implies the result.
\end{proof}

\begin{proof}[Proof of Lemma~\ref{small}.]
  This lemma is standard for Lelong numbers (see \cite{De2}). 
  We emphasize that the inequality (2) is not obvious and relies 
  on the Ohsawa-Takegoshi extension theorem. 
  
  The generalization is straightforward for Kiselman numbers.
  We nevertheless give the arguments for assertion (1) for completeness.
  
  Let $h \in \cH (s u)$ normalized by $| h|_{su} = 1$. The mean value property
  inequality for subharmonic functions implies
  \begin{multline*}
    |h(z,w)|^2
    \le\frac{C_1}{r^{2(1/\a_1+1/\a_2)}}
    \int_{|z-\xi|<r^{1/\a_1},~|w- \zeta|< r^{1/\a_2}}
    | h (\xi,\zeta)|^2 d\xi d\zeta
    \\
    \le \frac{C_1^{1/2}}{r^{2(1/\a_1+1/\a_2)}}
    \sup_{|z-\xi| < r^{1/\a_1},~|w- \zeta|< r^{1/\a_2}}e^{2 s u}.
  \end{multline*}
  Hence for all $(z,w) \in B$, we have
  \begin{equation*}
    u_s (z,w) 
    \le  
    s^{-1} \log \left(  \frac{C_1}{r^{1/\a_1+ 1/\a_2}}  \right) 
    + \sup _{|z-\xi| < r^{1/\a_1},~|w- \zeta|< r^{1/\a_2}} u.
  \end{equation*}
  Write $p =(z_0, w_0)$. By definition of Kiselman number we have
  \begin{equation*}
    \sup_ {|z_0-\xi| < r^{1/\a_1},~|w_0- \zeta|< r^{1/\a_2}}
    u \le (\a_1 \a_2)^{-1}~\nu ( p, u, (\a_1, \a_2))~\log r + C_3
  \end{equation*}
  hence
  \begin{multline*}
    \sup_ {|z_0-\xi| < r^{1/\a_1},~|w_0- \zeta|< r^{1/\a_2}} u_s
    \le 
    \sup_ {|z_0-\xi| < 2r^{1/\a_1},~|w_0- \zeta|< 2r^{1/\a_2}} u
    \\
    \le  
    (\a_1 \a_2)^{-1} \log (Cr )~\nu (p, u, (\a_1, \a_2) ) + C_3 - s^{-1}
    \log\left(\frac{C_1}{r^{1/\a_1+1/\a_2}}\right)
  \end{multline*}
  with $C = 2^{- 1/ \max \{ \a_1, \a_2 \} }$.
  We hence get
  \begin{equation*}
    \nu (p, u_s, (\a_1, \a_2) ) \ge \nu (p, u, (\a_1, \a_2) ) -
    s^{-1} (\a_1 + \a_2 )
  \end{equation*}
  as desired.
\end{proof}

%
%
%
%
\section{Asymptotic Multiplicities}\label{multi}

As explained earlier, the main ingredient in the proof of
Theorem~A are estimates from below of $\vol f^nB$ for a ball $B$.
This is equivalent to estimates of the Jacobian $Jf^n$, and more
precisely to the asymptotic order of vanishing of $Jf^n$. In this
section, we also consider two other
asymptotic multiplicities for a holomorphic mapping of $\PP{2}$.

We will make frequent use of the following 
\emph{strong Birkhoff theorem} due to the first author 
(see~\cite{Fav3} or~\cite{Fav4} Theorem 2.5.14):
\begin{Thm}\label{MT1}
  Let $f:\PP{2}\self$ be holomorphic of degree $d\ge2$ and
  let $\kappa_n:\PP{2}\to [ 1, + \infty [ $ 
  be a sequence of functions satisfying the following conditions:
  \begin{enumerate}
  \item
    for any $n \ge 0$, $\kappa_n$ is upper semicontinuous (usc) 
    with respect to the analytic Zariski topology;
  \item
    For any $n,m \ge 0$ and any $ x \in \PP{2}$, 
    \begin{equation*}    
      \kappa_{n+m}(x)\le\kappa_n(x)\,\kappa_m(f^nx);
    \end{equation*}
    we say $\kappa_n$ defines a submultiplicative cocycle;
  \item 
    for any $n \ge 0$, $\min_{\PP{2}} \kappa_n =1$.
  \end{enumerate}
  Then, for any $x \in \PP{2}$, the sequence $\kappa_n (x)^{1/n}$
  converges. Let  $\kappa_{\infty}(x)$ be its limit. We have
  $\kappa_{\infty} \circ f = \kappa_{\infty}$.
  Further, if $\kappa_{\infty}(x) > 1$, then
  \begin{itemize}
  \item
    either $x$ is preperiodic,
  \item
    or some iterate of $x$ belongs to a (not necessarily irreducible) 
    fixed curve $V$
    such that $\min_V\kappa_{\infty}=\kappa_{\infty}(x)$.
  \end{itemize}
\end{Thm}

\begin{Remark}
Note that the curve $V$ of the preceding theorem must contain an
irreducible component of  
the proper analytic subset $\cC \= \{ \kappa_1 >1 \}$.
\end{Remark}

%
%
\subsection{Asymptotic multiplicity of the Jacobian}

First we study the asymptotic order of vanishing of the Jacobian
of $f$.
Fix local charts $U\ni p$, $V\ni fp$, and denote by $Jf$ 
the Jacobian determinant of $f:U \to V$. 

\begin{Def}\label{asy:d-1}
  Let $\mu (p , Jf) \in \N$ be the order of vanishing 
  of $Jf$ at $p$.
\end{Def}
This number does not depend on the choice of chart. It can be
interpreted analytically as the Lelong number
\begin{equation}\label{asy:e-1}
  \mu(p,Jf)=\nu(p,dd^c\log|J f|)
\end{equation}
of the positive closed $(1,1)$ current $dd^c\log |Jf|$.
Note that $\mu(p,Jf)\ge 1$ if and only if $p$ belongs to the critical set
$\cC_f$.

We are interested in  studying the 
growth of $\mu(p,Jf^n)$ when $n$ tends to
infinity. It is straightforward to see that
\begin{equation}\label{asy-e:2}
  \mu(p,Jf^{k+n})=\mu (p,Jf^n)+\mu(p,Jf^k\circ f^n)
\end{equation}
for any $n,k\ge0$. The sequence $\{ \mu(p,Jf^n) \}_{n\ge 0}$ is not 
submultiplicative, but we will see that we can treat it as such.
First we have the following inequality.
\begin{Prop}{(see~\cite{Fav2} Remark~3.)}\label{asy-p:1}

  For any $p \in \PP{2}$ and any $n,k \ge 0$, the following inequality
  \begin{equation}\label{asy-e:3}
    \mu(p,Jf^{k}\circ f^n)\le(3+2\mu(p,Jf^n))\cdot\mu(f^np,Jf^k)
  \end{equation}
holds.
\end{Prop}
From this it follows that
\begin{align*}
  3+2\mu(p,Jf^{k+n}) 
  &=3+2\mu(p,Jf^n)+2\mu(p,Jf^k\circ f^n)\\
  &\le3+2\mu(p,Jf^n)+(3+2\mu(p,Jf^n))~(2\mu(f^np,Jf^k))\\
  &\le(3+2\mu(p,Jf^n))\cdot(3+2\mu(f^np,Jf ^k)).
\end{align*}
Introduce
\begin{equation*}
  \wmu(p,Jf^n)\=3+2\mu(p,Jf^n).
\end{equation*}
The last inequality can then be rewritten as
\begin{equation}\label{asy-e:5}
  \wmu(p,Jf^{k+n})\le\wmu(p,Jf^n)\cdot\wmu(f^np,Jf^k).
\end{equation}
The sequence $\wmu (p, Jf^n)$ hence defines an submultiplicative
cocycle. It is moreover usc with respect to the
analytic Zariski topology on $\PP{2}$ (e.g.\ by Siu's theorem). 
Thus Theorem~\ref{MT1} applies and yields:
\begin{Prop}\label{asy-p:2}
  Let $f:\PP{2}\self$ be a holomorphic map of degree $d\ge 2$.
  For any $p\in\PP{2}$, the sequence $\wmu(p,Jf^n)^{1/n}$ converges to
  a real number $\mu_{\infty}(p)\ge 1$. 
  We have $\mu_{\infty}\circ f=\mu_{\infty}$.
  Further, if $\mu_{\infty}(p)>1$, then one of the following holds:
  \begin{itemize}
  \item[(i)]
    $f^Np$ is a periodic critical point for some $N\ge0$;
  \item[(ii)]
    there exists a fixed curve $V$
    such that $f^Np\in V$ for some $N\ge 0$, and 
    $\min_V\mu_{\infty}=\mu_{\infty}(p)$.
  \end{itemize}
\end{Prop}

\begin{Remark}
  The sequence $ (1+\mu(p, Jf^n))^{1/n}$ also converges to $\mu_{\infty}(p)$.
  Moreover, define $\mu_n(p):=\mu(p,Jf\circ f^n)$. This sequence is clearly
  increasing. From~\eqref{asy-e:2}, it follows that
  \begin{equation*}
    n^{-1} \mu(p,Jf^n)\le\mu_{n+1}(p)\le\mu(p,Jf^n)~,
  \end{equation*}
  and so $(1+ \mu_n(p))^{1/n}$ also converges towards $\mu_{\infty}(p)$.
\end{Remark}

\begin{Prop}\label{asy-p:5}
  Let $f:\PP{2}\self$ be holomorphic of degree $d\ge2$. Then
  \begin{equation}\label{asy-e:4}   
    0\le\mu(p,Jf)\le 3(d-1)
  \end{equation} 
  for all $p\in\PP{2}$.
  Hence
  \begin{equation}   \label{asy-e:6}   
    1\le\mu_{\infty}(p)\le d
  \end{equation} 
  for all $p\in\PP{2}$.
\end{Prop}

\begin{proof}
  The multiplicity $\mu(p,Jf)$ is always smaller than the
  degree of the critical set of $f$, which is $3(d-1)$.
  Applying~\eqref{asy-e:4} to $f^n$ and letting $n \to \infty$
  we get~\eqref{asy-e:6}.
\end{proof}

%
%
\subsection{Asymptotic topological degree}

An important quantity in the stu\-dy of totally invariant sets is the local 
topological degree $e(p,f)$ of $f$ at a point $p$. 
By definition, this is the topological degree of the germ of an open map
induced by $f$ at $p$. 
Clearly $e(p,f)>1$ if and only if $p\in\cC_f$ and 
moreover $e(p,f)$ is an usc 
function with respect to the Zariski topology (see e.g.~\cite{Fav3}). Further
$e$ satisfies the composition formula
\begin{equation*}
  e(p,f^{k+n})=e(p,f^n)\cdot e(f^np,f^k).
\end{equation*}
Theorem~\ref{MT1} again applies and shows that
\begin{Prop}\label{asy-p:3}
  Let $f:\PP{2}\self$ be a holomorphic map of degree $d\ge 2$.
  For any $p\in\PP{2}$, the sequence $e(p,f^n)^{1/n}$ converges to
  a real number $e_{\infty}(p)\ge 1$. 
  We have $e_{\infty}\circ f=e_{\infty}$.
  Further, if $e_{\infty}(p)>1$, then one of the following holds:
  \begin{itemize}
  \item[(i)]
    $f^Np$ is a periodic critical point for some $N\ge0$;
  \item[(ii)]
    there exists a fixed curve $V$
    such that $f^Np\in V$ for some $N\ge 0$ and $\min_V e_{\infty} =
    e_{\infty}(p)$.
  \end{itemize}
\end{Prop}

\begin{Prop}\label{MP2}
  Let $f:\PP{2}\self$ be holomorphic of degree $d\ge2$. Then
  \begin{equation}\label{asy-e:7}   
    1\le e(p,f)\le d^2
  \end{equation}
  for all $p\in\PP{2}$ and $\{e(p,f)>d\}$ is a finite set whose 
  cardinality is bounded only in terms of the degree $d$.
  We have 
  \begin{equation}\label{asy-e:8}   
    1\le e_{\infty}(p)\le d^2
  \end{equation}
  for all $p\in\PP{2}$ and the set
  $\{e_\infty(p)=d^2\}$ is finite, totally invariant and
  contained in $\{e(p,f)=d^2\}$.
\end{Prop}

\begin{proof}[Proof of Proposition~\ref{MP2}.]
  Since the (global) topological degree of $f$ is $d^2$ we have
  $1\le e(p,f)\le d^2$ for all $p$. 
  Let us show that $ \{ e> d \}$ is a finite set.
  We follow the proof of Theorem~4.7 in~\cite{FoSi1}. 
  Take a point $z$ is the regular
  locus of $\cC_f$ such that $fz$ belongs to the regular locus of
  $f\cC_f$, and $z$ is a regular point for $f|_{\cC_f}$. 
  One can find local coordinates so that 
  $f(z,w)=(z,w^e)$ with $e \= e(p,f)$.
  On the other hand, we have in terms of
  currents 
  $f_*  [\cC_f]  = e [f\cC_f]$ and $ f^* f_* [\cC_f] =  e f^* [f\cC_f]
  \ge e^2 [\cC_f]$.
  As $\deg ( f_* [\cC_f]) = d \deg (\cC_f)$, we get $e^2 \le d^2$.
  Hence outside the finite set 
  $E\=\sing(\cC_f)\cup f^{-1}\sing(f\cC_f)\cup\sing(f|_{\cC_f})$ 
  we have $e\le d$.
  We conclude the proof noting that the cardinality of $E$ can be
  bounded only in terms of $d$.
  
  If $e_\infty(p)=d^2$, then the orbit of $p$ must visit the 
  finite set $\{e=d^2\}$ infinitely many times and is therefore
  preperiodic to a periodic orbit in that set. So we may write
  $f^mp=q=f^nq$ with $e(q,f)=d^2$. But then 
  $e_\infty(p)=e(q,f^n)^{1/n}$, so we must have 
  $e(f^iq,f)=d^2$ for $0\le i<n$. Thus 
  $f^{-1}(f^{i+1}q)=\{f^{i}q\}$ for all $i$, and that implies
  that $p=f^iq$ for some $i$. We conclude that $p$ 
  is periodic and that the orbit of $p$ is totally invariant
  and contained in $\{e=d^2\}$.
\end{proof}

%
%
\subsection{Asymptotic diameter}

A third quantity that we will need is related to the asymptotic 
diameter of balls $f^nB$.

To define this, fix local coordinates around $p$ and $fp$ so that
$p=fp=0$. Define $c(p,f)$ to be the largest integer $c$ such that
\begin{equation*}
  \vert f(\zeta)\vert\le A\vert\zeta\vert^c\qquad\text{as $\zeta\to0$}
\end{equation*}
for some constant $A>0$. Alternatively, $c(p,f)$ is the order of
the first nonvanishing term in the Taylor expansion of $f$. 
It can also be interpreted as the Lelong number $c(p,f) = \nu ( p,
\log |f|)$; hence the
definition does not depend on the choice of local coordinates and $c$
is usc with respect to the Zariski topology.
It is clear from the definition that
\begin{equation}\label{asy-e:9}
  c(p,f^{k+n})\ge c(p,f^n)\cdot c(f^np,f^k),
\end{equation}
and so $c$ defines a supermultiplicative cocycle. We have
\begin{Lemma}\label{asy-l:1}
  $\mu(p,f)\ge 2(c(p,f)-1)$.
\end{Lemma}
\begin{proof}
  This is a local result so we may assume $p=fp=0$.
  Write  $c=c(p,f)$ and $\mu=\mu(p,f)$.
  Let $f(\zeta)=f_c(\zeta)+O(|\zeta|^{c+1})$ where $f_c$ is a 
  homogeneous polynomial of degree $c$. The Jacobian determinant of
  $f_c$ is a homogeneous polynomial of degree $2c-2$ or vanishes identically.
  Thus $Jf(\zeta)=Jf_c(\zeta)+O(|\zeta|^{2c-1})$ and $\mu\ge2c-2$.
\end{proof}
This estimate and~\eqref{asy-e:9} 
allow us to deduce the following result from Proposition~\ref{asy-p:2}:
\begin{Prop}\label{asy-p:4}
  Let $f:\PP{2}\self$ be a holomorphic map of degree $d\ge 2$.
  For any $p\in\PP{2}$, the sequence $c(p,f^n)^{1/n}$ converges to
  a real number $c_{\infty}(p)\ge 1$. 
  We have $c_{\infty}\circ f=c_{\infty}$.
  Further, if $c_{\infty}(p)>1$, then one of the following holds:
  \begin{itemize}
  \item[(i)]
    $f^Np$ is a periodic critical point for some $N\ge0$;
  \item[(ii)]
    there exists a fixed curve $V$
    such that $f^Np\in V$ for some $N\ge 0$, 
    and $\min_Vc_{\infty}=c_{\infty}(p)$.
  \end{itemize}
\end{Prop}

\begin{Prop}\label{MP3}
  Let $f:\PP{2}\self$ be holomorphic of degree $d\ge2$. Then
  \begin{align}    
    1 & \le  c(p,f)  \le d, \label{asy-e:10}\\
    1 & \le  c_{\infty}(p,f)  \le d. \label{asy-e:12}
  \end{align}
  for any $p \in \PP{2}$.
  Moreover $c(p,f)=d$ if and only if
  $p$ is a homogeneous point for $f$, i.e. $f$ maps the 
  pencil of lines through $p$ to the pencil of lines through $fp$.
\end{Prop}

\begin{proof}
  By pre- and post-composing by projective linear maps of
  $\PP{2}$ we may assume that $p=fp=0$. Write
  $f(\zeta)=(P(\zeta)/R(\zeta),Q(\zeta)/R(\zeta))$
  where $P$, $Q$ and $R$ are polynomials of degree $d$ and
  $P(0)=Q(0)=0$, $R(0)=1$. Then clearly $f$ can only vanish
  up to order $d$ and so $1\le c(p,f)\le d$. Further, 
  $c(p,f)=d$ if and only if $P$ and $Q$ are homogeneous polynomials
  of degree $d$, which means precisely that $p$ is a homogeneous
  point.
\end{proof}

%
%
\subsection{Properties of the multiplicities}

We summarize in the following proposition the inequalities relating
the multiplicities $\mu, c$, and $e$ considered above.

\begin{Prop}\label{asy-p:7}
  Let $f:\PP{2}\self$ be holomorphic of degree $d\ge2$. For all
  $p\in\PP{2}$, we have
  \begin{alignat}{1}
    2(c(p,f)-1) &\le  \mu(p,Jf) \le 2~(e(p,f)-1),  
    \label{ineq1}
    \\ 
    c(p,f)  &\le  \sqrt{e(p,f)}.
    \label{ineq2}
  \end{alignat}
  Hence 
  \begin{alignat}{2}
    c_\infty(p) \le \mu_{\infty}(p) &\le    e_\infty(p) &\le d, 
    \label{ineq3}
    \\
    c_\infty(p) &\le \sqrt{e_\infty(p)} &\le d.
    \label{ineq4}
  \end{alignat}
  The set $\{c_\infty=d\}\subset\{e_\infty=d^2\}$ is finite
  and totally invariant.
\end{Prop}

\begin{Example}\label{ME1}
  A totally invariant point, i.e.\ a point with $e(p,f)=d^2$,
  is not necessarily superattracting as the
  example $f(z,w)= (2z+w^d,z^d)$ from~\cite{FoSi1} shows.
  In this example the origin is a totally 
  invariant fixed point with one expanding eigenvalue $2>1$.
  One can check that $c_\infty=\mu_\infty=1$, whereas
  $e_\infty=d^2$ for this map.
\end{Example}

\begin{proof}[Proof of Proposition~\ref{asy-p:7}.]
  Equations~\eqref{ineq3} and~\eqref{ineq4} are consequences
  of~\eqref{ineq1} and~\eqref{ineq2}, 
  and the last assertion follows from~\eqref{ineq4} 
  and Proposition~\ref{MP2}.
  
  Equations~\eqref{ineq1} and~\eqref{ineq2} 
  are local so we may assume $p=fp=0$ and for sake
  of simplicity we write
  $\mu(p,f) = \mu$, $c(p,f) =c$, $e(p,f) =e$.
  Because of Lemma~\ref{asy-l:1} we only have to show the inequalities
  $\mu \le 2 e -2$ and $ c \le \sqrt{e}$.
  
  To prove 
  $\mu \le 2(e-1)$ we use the fact that
  $|f(\zeta)|\ge C|\zeta|^e$ (see e.g.~\cite{Fav2}) and (by definition)
  $|Jf(\zeta)|\le D|\zeta|^{\mu}$ for some constants $C,D >0$.
  For any ball $B(r)$ of radius $r>0$, one gets
  \begin{align*}
    \vol fB(r) &\ge \vol B(Cr^e)=C'r^{4e}\\
    \vol fB(r) &= e^{-1}\int_{B(r)}|Jf|^2 
    \le e^{-1}\int_{B(r)}D|\zeta|^{2\mu}=D'r^{2\mu+4} 
  \end{align*}
  for some constants $C',D'>0$. By letting $r\to0$ we get
  $\mu\le2(e-1)$. 
  
  We now give an analytic proof of $c \le  \sqrt{e}$. We first note that 
  $\delta_p = ( dd^c \log |\zeta|)^2$.
  It follows from~\cite{De1} Corollary 6.8 that 
  \begin{align*}
    e 
    & =   
    \nu ( f^* \delta_p , p) 
    \\
    &  =   
    \nu ( f^* dd^c \log |\zeta| \wedge f^* dd^c \log |\zeta| , p)  
    \\
    & \ge  
    \nu (f^* dd^c \log |\zeta|, p )^2
    = c^2.
  \end{align*}
  This concludes the proof.
\end{proof}

%
%
\subsection{Exceptional sets}

We define two exceptional sets as follows.
\begin{Def}\label{asy-d:2}
  Let $f:\PP{2}\self$ be holomorphic of degree $d$.
  \begin{itemize}
  \item
    Define \emph{the first exceptional set} 
    $\cE_1$ to be the union of irreducible 
    curves $V$ such that $\mu_{\infty}(p)=d$ for all $p\in V$.
  \item
    Define \emph{the second exceptional set} $\cE_2$ 
    to be the set of points $p$ with $c_\infty(p)=d$.
  \end{itemize}
  Finally, we define the \emph{exceptional set} $\cE$
  by $\cE:=\cE_1\cup\cE_2$.
\end{Def}

The exceptional set $\cE$ is where $f$ is most volume contracting; 
we will spend the next few sections analyzing it.
In particular we will show that $\cE$ is algebraic, totally
invariant and superattracting.
%
%
%
%
\section{The first exceptional set $\mathcal{E}_1$}\label{exc1}

The key to the description of the first and second
exceptional sets $\cE_1$ and $\cE_2$ lies in understanding 
the loci $\{c_\infty<\mu_\infty=d\}$ and $\{c_\infty=d\}$.
Being in the first locus means, roughly speaking, that volume is 
contracted much faster than diameter. This implies that the image of a
ball is very close to being one-dimensional, and that can only 
happen if the critical set is very recurrent. Being 
in the second locus means that diameter is decreasing very fast,
and this leads to totally invariant and superattracting points.
In fact, we already know from Proposition~\ref{asy-p:7} that the 
second exceptional set 
$\cE_2=\{c_\infty=d\}$ is finite and totally invariant.
In this section, we show

\begin{Thm}\label{str-t:2}
 Let $f:\PP{2}\self$ be a holomorphic map of degree $d\ge2$.
 Then the first exceptional set $\cE_1$ consists of the union of (not
 necessarily irreducible) totally invariant curves,
 and equals the (Zariski)
 closure of the locus $\{c_\infty<\mu_\infty=d\}$.
 Moreover, $\cE_1$ is the union of at most three lines 
 in general position. 
\end{Thm} 

The proof relies essentially on the following local result.
Notice that $c_\infty$ and $\mu_\infty$ can be defined for a  germ
fixing a point.
\begin{Thm}\label{str-t:1}
  Let $f:(\CC{2},0)\self$ be a holomorphic germ. 
  Let $V_1,\dots,V_k$ be the irreducible components of the 
  critical set $\cC_f$. Assume that 
  $c_{\infty}(0)<\mu_{\infty}(0)$.
  Then there exist $a_1,\dots,a_k\ge0$ such that
  \begin{equation}\label{str-e:1}
    f^*\left(\sum_ia_i[V_i]\right)
    \ge 
    \mu_{\infty}(0)\left(\sum_i a_i[V_i]\right).
  \end{equation}
\end{Thm}

\begin{proof}[Proof of Theorem~\ref{str-t:2}.]
  First assume that $V$ is a totally invariant curve.
  We want to show $V\subset\cE_1$.
  Given $n\ge1$ and $p\in V$ outside a finite subset 
  (depending on $n$) we may pick local coordinates at $p$ and
  at $f^np$ so that $f^n$ is given by $(z,w)\mapsto(z^{d^n},w)$.
  It follows that $\mu(p,Jf^n)=d^n-1$ on $V$ outside a finite set, 
  and thus $\mu(p,Jf^n)\ge d^n-1$ on all of $V$ by upper semicontinuity.
  This implies $\mu_\infty=d$ on $V$ hence $V\subset\cE_1$.

  Conversely pick an irreducible component $V\subset\cE_1$. We will
  show that $V$ belongs to a totally invariant curve. 
  We claim that
  \begin{equation}\label{cla}
    \{ c_{\infty}<  \mu_{\infty} = d  \} \subset \cE_1.
  \end{equation}
  This implies that $V\setminus\{c_{\infty}=d\} \subset \cE_1$,
  hence $V \subset \cE_1$ as $\{c_{\infty}=d\}$ is finite.
  
  Let us show the claim.
  Pick a point $p$ satisfying $c_{\infty}(p)<\mu_{\infty}(p)=d$. 
  By Proposition~\ref{asy-p:2}, either
  $p$ is preperiodic or $f^mp \in W$ for some $ m \ge 0$ and for some
  fixed curve $W$ with $\min_W \mu_{\infty} = \mu_{\infty}(p) = d$.  
  
  In the latter case, we have $f^* [W] \ge l [W]$ for some 
  maximal $l\in[2,d]$. For a generic point $q\in W$, we have $d =
  \mu_{\infty}(q) = e_{\infty}(q)$, and for any $j \ge 0$, $e ( q , f^j) =
  l^j$. Hence $l =d$, and
  we infer that $W$ is totally invariant, whence $p \in W \subset \cE_1$.
  
  If $p$ is preperiodic, then
  $f^mp=q=f^Nq$ for some $q\in\cC_f$ and $m\ge0$. 
  Clearly $\mu_\infty(q,f^N)=\mu_\infty(p,f)^N=d^N$ and
  $c_\infty(q,f^N)=c_\infty(p,f)^N<d^N$.
  Apply Theorem~\ref{str-t:1} to $f^N$ and
  find non-negative integers $a_1,\dots,a_k$ such that
  $f^{N*}\sum_ia_i[V_i]\ge d^N\sum_ia_i[V_i]$.
  As $f^N$ is of degree $d^N$, we have equality. In particular 
  the union $W$ of the critical components passing through $p$ with 
  $a_i>0$ is a totally invariant set for $f^N$. But then
  the curve $W'=W\cup\dots\cup f^{N-1}W$ is
  totally invariant for $f$, hence $W' \subset \cE_1$ by the argument
  above. Since $f^mp\in W'$ we have in fact
  $p\in W'$. 
  This concludes the proof of the claim.
  
  Finally,~\eqref{cla} shows that
  $\cE_1$ is the closure of the set
  $\{ c_{\infty} < \mu_{\infty} = d \}$.
  The remaining statement of the theorem follows from the classification
  of totally invariant curves (Proposition~\ref{MP1}).
\end{proof}

\begin{proof}[Proof of Theorem~\ref{str-t:1}.]
  Pick holomorphic maps $\phi_i$ so that  
  $V_i=\phi_i^{-1}(0)$, and set $\phi:=Jf$. Then
  $\phi_i$ are the irreducible factors of $\phi$.
  There exist integers $t_{ij}\ge 0$ so that
  \begin{equation}\label{str-e:2}
    \phi_i\circ f=\wphi_i\times\prod_{j=1}^k\phi_j^{t_{ij}} 
  \end{equation}
  for some holomorphic $\wphi_i$ with 
  $\wphi_i^{-1}(0)\cap\cC_f=\{0\}$.
    By {\L}ojasiewic's inequality~\cite[p. 243]{Lo1}
  there exist constants $C,\alpha>0$ such that
  \begin{equation}\label{str-e:3}
    |\wphi _i(\zeta)|+|\phi(\zeta)|\ge C|\zeta|^{\alpha}
  \end{equation}
  in a neighborhood of the origin.
  For fixed $n\ge 0$ write 
  $f^n=f_{c_n}+ O(|\zeta|^{c_n +1})$
  where $f_{c_n}\not\equiv0$ 
  is a homogeneous polynomial of degree $c_n$.
  Similarly, set 
  $\phi\circ f^n=\phi_{\mu_n}+O(|\zeta|^{\mu_n+1})$ 
  for a non-degenerate homogeneous polynomial $\phi_{\mu_n}$
  of degree $\mu_n=\mu(0,Jf \circ f^n)$.
  
  We know that $c_n$ is 
  an increasing supermultiplicative sequence such that 
  $c_n^{1/n}\to c_{\infty}(0)$, and that $\mu_n$ is an
  increasing sequence such that $\mu_n ^{1/n} \to \mu_{\infty}(0)$ 
  (see Section~\ref{multi}). 
  By hypothesis $c_{\infty}(0)<\mu_{\infty}(0)$, so for any 
  fixed $c$, $\mu$ with
  $c_{\infty}(0)<c<\mu<\mu_{\infty}(0)$ we have
  \begin{align*}
    c_n &\le A c^n \\
    \mu_n &\ge B\mu^n 
  \end{align*}
  for $n\ge 0$ and some constants $A,B>0$.
  
  We infer that for a generic $|\zeta|\ll1$ and for any $n\ge 0$
  \begin{equation}\label{str-e:4}
    |\zeta_n|:=|f^n(\zeta)|\ge C_n|\zeta|^{Ac^n}
  \end{equation}
  for some $C_n >0$.
  On the other hand for all $|\zeta|\ll1$ and all $n\ge 0$ we have
  \begin{equation}\label{str-e:5}
    |\phi\circ f^n(\zeta)|=|\phi(\zeta_n)|\le D_n|\zeta|^{B\mu^n}.
  \end{equation}
  
  Let $M$ be the $k$ by $k$ matrix $[t_{ij}]$ and let $\rho>0$ be the
  spectral radius of $M$. Denote $M^n \= [t_{ij}^n]$ and fix $K>0$ so that
  \begin{equation}\label{stru-e:6}
    0\le t_{ij}^n \le K\rho^n.
  \end{equation}
  For $n$ large enough, so that $Ac^n\le B\mu^n$, and for generic
  $\zeta$, we can apply~\eqref{str-e:3}, and get
  \begin{equation*}
    |\wphi_j(\zeta_n)|
    \ge 
    C|\zeta_n|^{\alpha}-|\phi(\zeta_n)|
    \ge
    CC_n^{\alpha}|\zeta|^{Ac^n}-D_n|\zeta|^{B\mu^n}
    \ge 
    C_n'|\zeta|^{A c^n}.
  \end{equation*}
  Hence 
  \begin{equation*}
    \left[
      \begin{array}{c}
        |\phi_1|
        \\
        \vdots
        \\
        |\phi_k|
      \end{array}
    \right]
    (\zeta_{n+1})
    =
    \left[
      \begin{array}{c}
        |\wphi_1|\prod_j|\phi_j|^{t_{1j}}
        \\
        \vdots
        \\
        |\wphi_k|\prod_j|\phi_j|^{t_{kj}}
      \end{array}
    \right]
    (\zeta_{n})
    \ge 
    C'_n |\zeta|^{A c^n} 
    \left[
      \begin{array}{c}
        |\phi_1|
        \\
        \vdots
        \\
        |\phi_k|
      \end{array}
    \right]^M
    (\zeta_{n}),
  \end{equation*}
  where we let 
  ${}^t[|\phi_1|,\dots,|\phi_k|]^M
  \={}^t[\prod_j|\phi|_j^{t_{1j}},\dots,\prod_j|\phi|_j^{t_{kj}}]$.
  
  By induction and using~\eqref{stru-e:6} for $1\le i\le k$ we have
  \begin{equation*}
    |\phi_i(\zeta_{n+1})|
    \ge C''_n |\zeta|^{A K \frac{\rho^{n+1}- c^{n+1}}{\rho - c}}.
  \end{equation*}
  On the other hand~\eqref{str-e:5} shows that
  \begin{equation*}
    D_{n+1}|\zeta|^{B\mu^{n+1}} 
    \ge|\phi(\zeta_{n+1})|
    \ge C'''_n|\zeta|^{A'K\frac{\rho^{n+1}-c^{n+1}}{\rho-c}}.
  \end{equation*}
  As $\zeta$ was generic, we can let it tend to zero and let $n$ tend to
  infinity. We infer $\rho \ge \mu$, and therefore $\rho \ge
  \mu_{\infty}(0)$ as $\mu < \mu_{\infty}(0)$ was chosen arbitrary.
  
  The Perron-Frobenius theorem now implies the existence of an
  eigenvector $(a_1,\dots,a_k)$ for $M$ with non-negative
  coefficients associated to the eigenvalue $\rho$.
  We have $ f^* \sum_i a_i [V_i] = \rho ( \sum_i a_i [V_i]) \ge
  \mu_{\infty}(0) \sum_i a_i [V_i]$, which completes the proof.
\end{proof}
                               
%
%
%
%
\section{The second exceptional set $\mathcal{E}_2$}\label{exc2}

The second exceptional set $\cE_2$ is both hard and interesting
to analyze in detail. We will give some partial results that
are enough for the purpose of Theorems~A and~A'.
\begin{Prop}\label{MP5}
  The second exceptional set is given by the set of 
  \begin{itemize}
  \item
    points $p\notin\cE_1$ with $\mu_{\infty}(p)=d$;
  \item
    totally invariant periodic orbits $f$ in $\cE_1$.
  \end{itemize}
  Further, $\cE_2$ is finite, totally invariant and 
  superattracting.
\end{Prop}
\begin{proof}
  Proposition~\ref{asy-p:7} shows that
  $\cE_2=\{c_\infty=d\}$ is finite and totally invariant.
  A point $p\in\cE_2$, $f^np=p$, is superattracting as
  the Taylor series of $f^n$ at $p$ 
  vanishes to order $c(p, f^n) \ge 2$. 
  Let us now prove the characterization of $\cE_2$.
  First consider a point $p\notin\cE_1$. If $p\in\cE_2$, then
  $c_\infty(p)=d$ and so $\mu_\infty(p)=d$ by 
  Proposition~\ref{asy-p:7}. Conversely, if $\mu_\infty(p)=d$,
  then $c_\infty(p)=d$ by Theorem~\ref{str-t:2} so
  $p\in\cE_2$.

  Next consider $p\in\cE_1$. We want to show that $c_\infty(p)=d$
  if and only if $p$ is periodic and the orbit of $p$ is totally
  invariant. Both of these properties are preserved under replacing 
  $f$ by an iterate, so we may assume that the line $(w=0)$ is
  totally invariant for $f$ and that $p$ is on this line. Notice that
  the restriction $R \= f|_{w=0}$ is a rational map of
  degree $d$ and that $c(p,f^n)\le c(p,R^n)=e(p,R^n)$.

  Assume first $c_{\infty}(p)=d$. Then $e_{\infty}(p,R)=d$ hence $p$
  belongs to a totally invariant orbit.

  Conversely, suppose that the orbit of $p$ is totally invariant.
  After replacing $f$ by $f^2$ we may assume that $p=(0,0)$ is a
  totally invariant fixed point for $R$. 
  Thus we may write 
  \begin{equation*}
    f(z,w)=\left(
      \left(z^d+wQ\right)(1+\eta),
      w^d(1+\eta)
    \right)
  \end{equation*}
  in local coordinates $\zeta=(z,w)$, for holomorphic $Q$, $\eta$ 
  with $\eta(0,0)=0$.
  One checks that $c(p,f^n)\ge d^{n-1}$, thus $c_\infty(p)=d$.
  This completes the proof.
\end{proof}

Let us give some examples of maps with
nonempty second exceptional set.
\begin{Example}\label{ME2}
  We already know from Proposition~\ref{MP5}
  that a totally invariant point on a totally
  invariant line is in $\cE_2$. This will happen if
  $f[z:w:t]=[z^d+wQ(z,w,t):w^d:R(z,w,t)]$, where $Q$ and
  $R$ are homogeneous of degree $d-1$ and $d$, respectively.
\end{Example}
\begin{Example}\label{ME3}
  A homogeneous point, i.e.\ a point $p$ with the property that the
  family of lines through $p$ is invariant under $f$, is a point
  in $\cE_2$. In fact, $c(p,f^n)=d^n$ and so $c_\infty(p)=d$.
  In homogeneous coordinates where $p=[0:0:1]$ we have
  $f[z:w:t]=[P(z,w):Q(z,w):R(z,w,t)]$.
\end{Example}
\begin{Example}\label{ME4}
  More generally, if $f$ preserves a linear pencil of curves, then 
  any base point $p$ of the pencil belongs to $\cE_2$ for this map.
  Holomorphic maps of $\PP{2}$ preserving a pencil of curves were 
  studied in~\cite{DaJo1}: it turns out that the mappings and 
  the base points are always of one of the types described in
  Examples~\ref{ME2} or~\ref{ME3}.
\end{Example}  
It is a very interesting problem whether all points in $\cE_2$ are of
the types described in Examples~\ref{ME2} or~\ref{ME3}. We postpone
its discussion to a later paper.  Assuming this is the case, we refer
to the last section of the paper for a list of different possible
configurations of the exceptional sets $\cE_1$ and $\cE_2$.

%
%
%
%
\section{Volume estimates outside the exceptional set}\label{volume0}

In this section we give a lower bound on the volumes $f^nE$
for sets $E$ that avoid a fixed neighborhood of the exceptional
set $\cE$ under iterations.

\begin{Thm}\label{vol-t:1}
  Let $f:\PP{2}\self$ be a holomorphic map of degree $d\ge2$.
  Fix an open neighborhood $\Omega\supset\cE$ of the exceptional set.
  Then there exist a constant $\lambda<d$ and constants 
  $C_1,C_2>0$ such that 
  \begin{equation}\label{vol-e:1}
    \vol f^nE\ge(C_1\vol E)^{C_2 \lambda^n}
  \end{equation}
  for any Borel set $E\subset\PP{2}$ and any integer $n\ge 0$ with
  $E,\dots,f^nE\subset\PP{2}\setminus\Omega$.
\end{Thm}

The key idea ingredient in the proof is the following upper bound on
the multiplicities of the Jacobian.

\begin{Prop}\label{vol-p:1}
  
  In the setting of Theorem~\ref{vol-t:1} 
  there exist $\rho<d$ and $C>0$ such that
  \begin{equation}\label{vol-e:2}
    \mu(x,Jf^n)\le C\rho^n
  \end{equation}
  for any $x\notin\cE$ and any integer $n\ge 0$.
\end{Prop}

We defer the proof of Proposition~\ref{vol-p:1} until the end of this 
section and show how to deduce Theorem~\ref{vol-t:1} from it.
The key to doing so is the following result, which connects
the multiplicity of the Jacobian to the volume of $fE$.

\begin{Prop}{(see~\cite{Fav4} Chapitre 4)}\label{vol-p:2}
  Let $f:\PP{2}\self$ be a holomorphic mapping and
  $K\subset\PP{2}$ be a compact set. Define
  \begin{equation*}
    \tau_f(K):=\max\{\mu(p, Jf),~p\in K\}.
  \end{equation*}
  Then for any $\e>0$, 
  there exists a constant $C_{\e}>0$ such that
  \begin{equation}\label{vol-e:6}
    \vol fE\ge C_{\e}(\vol E)^{1+\tau_f(K)+\e} 
  \end{equation}
  for any Borel set $E \subset K$.
\end{Prop}
\begin{proof}
  Write $\tau_{\e}:=\tau_f(K)+\e$. 
  Let $E\subset K$ be a Borel set. 
  We are looking for a lower bound for
  $\vol fE $ in terms of $\vol E$. To this end we apply the
  Kiselman-Skoda estimate (Theorem~\ref{vol-t:2}) to the function 
  $\log|Jf|$ in each chart of a given atlas of $\PP{2}$.
  Notice that $\nu(z,\log|Jf|)=\mu(z,Jf)$.
  We conclude that there exists a constant $C_{\e}>0$ such that
  \begin{equation}\label{vol-e:7}
    \vol\left(K\cap\{|Jf|^2\le t\}\right)\le C_{\e}t^{1/\tau_{\e}}
  \end{equation}
  for all $t\ge0$.
  
  We pick a ``stopping time'' $T_0$ defined by
  $T_0^{1/\tau_{\e}}=(2C_{\e})^{-1}(1+\tau_{\e}^{-1})\vol E$, 
  and deduce the following sequence of inequalities:
  \begin{align*}
    \vol fE
    &\ge d^{-2}\int_{E}|Jf|^2 
    \ge d^{-2}\int_0^{T_0}\left(\vol E-\vol\{|Jf |^2\le t\}\right)dt\\
    &\ge d^{-2}\left(T_0 \vol E-\int_0^{T_0}C_{\e}t^{1/\tau_{\e}}\,dt\right)
    \ge C'_{\e}(\vol E)^{1+\tau_{\e}},
  \end{align*}
  which complete the proof.
\end{proof}

\begin{proof}[Proof of Theorem~\ref{vol-t:1}]
  Choose an integer $N$ so that $C\rho^N<d^N$ for the
  constant $\rho<d$ given by Proposition~\ref{vol-p:1}. 
  Fix $\l<d$ so that $C\rho^N<\lambda^N<d^N$.
  Proposition~\ref{vol-p:2} applied to $f^N$ with 
  $K\=\PP{2}\setminus\Omega$ yields a constant $C>0$
  so that
  \begin{equation}\label{vol-e:8}
    \vol f^NE\ge C(\vol E)^{\lambda^N}
  \end{equation} 
  for any Borel set $E\subset\PP{2}$.
  In a same way, one can find constants $D,D'>0$ so that
  \begin{equation}\label{vol-e:9}
    \vol f^jE\ge D'(\vol E)^{D}
  \end{equation}
  for any $0\le j\le N-1$ and any $E\subset\PP{2}$.
  
  Take a Borel set $E\subset\PP{2}$ and $n\ge 0$ so that 
  $E,\dots,f^nE\subset\PP{2}\setminus\Omega$.
  Write $n=kN+l$ with $l\ge N-1$. We have
  \begin{align*}
    \vol f^nE
    &= \vol f^l(f^{kN}E)\ge D'\vol f^{kN}E^{D}\\
    &\ge D'\left((\vol E)^{\lambda^{kN}}C^{\sum_{l=1}^{k-1}
        \lambda^{lN}}\right)^{D}
    \ge(C''\vol E)^{D''\lambda^n},
  \end{align*}
  which completes the proof of Theorem~\ref{vol-t:1}.
\end{proof}

Finally we prove the estimate for multiplicities 
in Proposition~\ref{vol-p:1}.
\begin{proof}[Proof of Proposition~\ref{vol-p:1}]
  Denote by $V_1,\dots,V_k$ the irreducible components of $\cC_f$
  that are not in $\cE_1$. For each $i$, pick a point $x_i \in V_i$ so
  that $\mu_{\infty}(x_i) <d$, and fix $\lambda < d$ with $\max _i
  \mu_{\infty}(x_i) < \l$. One can find a constant $C>0$ so that
  \begin{equation}\label{vol-e:3}
    \mu(x_i , Jf^n) \le C \lambda ^n 
  \end{equation}
  for all $n,i \ge0$.
  Introduce the set 
  ${\mathcal F}_N \= \{ x \in \PP{2} \setminus \cE,~\mu(x,Jf^N)>C\lambda^N \}$
  for a suitable $N \ge 0$ to be chosen later.
  Because of~\eqref{vol-e:3}, it is a finite set.
  
  Let $x \in \PP{2}$ and $n \ge 0$.
  
  $\bullet$ First assume 
  $\{ x,fx,\dots,f^nx\} \cap {\mathcal F}_N=\emptyset$. 
  Write $n = kN + l$ with $ 0 \le l \le N-1$.
  We have
  \begin{align*}
    \mu(x, Jf^n)
    &=\mu(x,Jf^{kN+l})=\mu(x,Jf^l)+\mu(x,Jf^{kN}\circ f^l)\\
    &\le\mu(x,Jf^l)+(3+2\mu(x,Jf^l))\mu(f^lx,Jf^{kN}).
  \end{align*}
  By~\eqref{asy-e:5}, we infer
  \begin{align*}
    3 + 2 \mu ( f^lx, Jf^{kN})
    \le 
    \prod_{j=0}^{k-1} 
    3 + 2 \mu ( f^{l+j}x, Jf^{N})
    \le ( 3 + 2 C \lambda^N) ^k
  \end{align*}
  Set $C_N \= \max _{\PP{2}} \mu(\cdot, Jf^N)$. Note that  $ C_N \ge
  C_{N-1}$.
  We get
  \begin{align*}
    \mu(x, Jf^n)
    & \le  
    C_N + ( 3 C_N + 2) \mu ( f^lx , Jf^{Nk})
    \\
    & \le  
    3^{-1} ( 3 C_N + 2)  ( 3 \mu ( f^lx , Jf^{Nk}) + 2 ) + C_N - 2/3 (
    3 C_N + 2) 
    \\
    & \le  
    (C_N + 2/3) ( 3 + 2 C \lambda ^N) ^k 
  \end{align*}
  Now for a fixed $\lambda < \rho < d$ we take $N \gg 1$ large enough 
  to conclude
  \begin{equation}\label{vol-e:5}
    \mu ( x, Jf ^n)
    \le 
    C' \rho ^n
  \end{equation}
  for some constant $C' >0$.
  
  $\bullet$ Now assume
  $\{x,fx,\dots,f^nx\}\cap{\mathcal F}_N\not=\emptyset$.
  
  The set ${\mathcal F}_N$ is finite. By definition it does not
  intersect $\cE$, hence one  can find constants $\lambda'<d$,
  $C''>0$ such that 
  \begin{equation}\label{vol-e:4}
    \mu (p , Jf^n) \le C'' (\l') ^n 
  \end{equation}
  for all $p \in {\mathcal F}_N$ and $n \ge 0$.

  Let $l$ be the smallest integer such that $ f^lx\in {\mathcal
    F}_N$. Applying~\eqref{asy-e:5} and~\eqref{vol-e:4}, we get
  \begin{align*}
    \mu(x, Jf^n)
    & \le   
    \mu ( x, Jf^l) 
    + 
    \mu ( f^lx, Jf^{n-l}) ( 3 + 2 \mu ( x, Jf^l))
    \\
    & \le  
    C' \rho ^l + ( 3 + 2 C' \rho ^l) ( C'' (\l') ^{n-l})
    \\
    & \le  
    D \max \{ \rho, \l' \} ^n,
  \end{align*} 
  and the proof is complete.
\end{proof}

%
%
%
%
\section{Volume estimates at the first exceptional set $\cE_1$}\label{volume1}

We now analyze the dynamics near $\cE_1$. Recall that $\cE_1$ is
a totally invariant union of at most three lines (in general position).

\begin{Prop}~\label{MP4}
  Let $f:\PP{2}\self$ be a holomorphic map of degree $d\ge2$.
  Fix small open neighborhoods $\Omega_i\supset\cE_i$ of the first and
  second exceptional sets, respectively. Then 
  there exists
  a constant $C>0$ and an integer $N\ge1$ such that 
  \begin{equation}\label{Me10}
    \vol f^nE\ge(C\,\vol E)^{d^n}
  \end{equation}
  for any Borel set $E\subset\PP{2}$ and any integer $n\ge N$ with
  $E,\dots,f^nE\subset\Omega_1\setminus\Omega_2$.
\end{Prop}

We will prove Proposition~\ref{MP4} using the structure of the
Jacobian $Jf$ near $\cE_1$. First we prove the following lemma.

\begin{Lemma}\label{ML2}
  There exists $C>0$ such that for all $s>0$ we have
  \begin{equation}\label{Me7}
    \vol\left\{p\in\Omega_1:\ |Jf(p)|<s\right\}
    \le Cs^\frac2{d-1}\log s.
  \end{equation}
  Further, there exists $N\ge 1$ and for any $n\ge N$ a constant
  $C_n>0$ such that 
  \begin{equation}\label{Me8}
    \vol\left\{p\in\Omega_1\setminus\Omega_2:\ |Jf^n(p)|<s\right\}
    \le C_ns^\frac2{d^n-1}
  \end{equation}
  for all $s>0$.
\end{Lemma}

\begin{proof}[Proof of Lemma~\ref{ML2}]
  The result is local. Pick $p\in\cE_1$. We may assume that $p=(0,0)$ and 
  that $\cE_1=(w=0)$ or $\cE_1=(zw=0)$ locally at $p$.
  Further, write 
  \begin{equation*}
    Jf(z,w)=w^{d-1}\prod_{i=1}^k(z-\alpha_i(w)),
  \end{equation*}
  where $1 \le k\le d-1$ is the multiplicity of $z=0$ as a 
  critical point of $f|_{(w=0)}$ and where $\alpha_i$ are
  multivalued functions with $\alpha_i(0)=0$. 
  Let $X_s$ be the subset of the bidisk 
  $\Delta^2=\{|z|,|w|<1\}$ where $|Jf|<s$. 
  We will show that 
  \begin{equation*}
    \vol X_s\le Cs^\frac2{d-1}\log s,
  \end{equation*}
  which will prove~\eqref{Me7}.
  For the rest of the proof we will let $C$ denote various
  positive constants.
  Let $A=\sup_{\Delta^2}|\prod_{i=1}^k(z-\alpha_i(w))|$. For 
  $\tau<A$ and fixed $w\in\Delta$ we may estimate
  \begin{align*}
    \area\left\{z\in\Delta:\ 
      \left|\prod_{i=1}^k(z-\alpha_i(w))\right|\le\tau\right\}
    &\le\sum_i\area\left\{\left|(z-\alpha_i(w))\right|\le\tau^{1/k}\right\}\\
    &\le C\tau^{2/k}.
  \end{align*}
  We can then use Fubini's Theorem to estimate the volume of $X_s$:
  \begin{align}\label{Me9}
    \vol X_s
    &\le 2\pi\int_0^{(s/A)^{\frac1{d-1}}}r\,dr\notag\\
    &+\int_{(s/A)^{\frac1{d-1}}}^1
    \area\left\{
      \left|\prod_{i=1}^k(z-\alpha_i(w))\right|<\frac{s}{r^{d-1}}
    \right\}r\,dr\notag\\
    &\le Cs^{\frac2{d-1}}
    +Cs^{\frac2{k}}\int_{(s/A)^{\frac1{d-1}}}^1
    r^{1-\frac{2(d-1)}{k}}\,dr.
  \end{align}
  If $k<d-1$, then the second term in~\eqref{Me9} can be estimated by
  \begin{equation*}
    Cs^{\frac2{k}}s^{\frac1{d-1}(2-\frac{2(d-1)}{k})}=Cs^{\frac2{d-1}},
  \end{equation*}
  and so $\vol X_s\le Cs^{2/(d-1)}$ in this case.
  
  If $k=d-1$ then the second term in~\eqref{Me9}
  is instead bounded by
  \begin{equation*}
    Cs^{\frac2{d-1}}\log(s^{\frac1{d-1}})=Cs^{\frac2{d-1}}\log s,
  \end{equation*}
  and so $\vol X_s\le Cs^{2/(d-1)}\log s$.
  This proves~\eqref{Me7}. 
  
  As for~\eqref{Me8} we notice that for $n\ge N$,
  all the critical points for $f^n|_{(w=0)}$, except the
  ones at $\cE_2$, will have multiplicity $<d^n-1$. Thus
  the above calculations imply~\eqref{Me8}.
\end{proof}

We now show how Lemma~\ref{ML2} implies Proposition~\ref{MP4}
\begin{proof}[Proof of Proposition~\ref{MP4}]
  The proof is similar to that of Proposition~\ref{vol-p:2}.
  First assume that $N\le n\le 2N$, with $N$ from Lemma~\ref{ML2}.
  We pick a ``stopping time'' $T_n$ defined by
  \begin{equation*}
    C_n(1+d^{-n})T_n^{\frac1{d^n-1}}=2^{-1}\vol E,
  \end{equation*}
  with $C_n$ from Lemma~\ref{ML2}.
  Then we get
  \begin{align*}
    \vol f^nE
    &\ge d^{-2n}\int_{E}|Jf^n|^2\\
    &\ge d^{-2n}\int_0^{T_n}
    \left(\vol E-\vol\{|Jf^n|^2<t\}\right)dt\\
    &\ge d^{-2n}
    \left(T_n\vol E-\int_0^{T_n}C_nt^{\frac1{d^n-1}}\,dt\right)\\
    &\ge   2^{-1}d^{-2n}T_n\vol E
    \ge C_n'(\vol E)^{d^n}.
  \end{align*}
  It is now easy to iterate this estimate and arrive at~\eqref{Me10}.
\end{proof}

%
%
%
%
\section{Attenuation of Lelong numbers}\label{lelong}

For the proof of Theorem~A we need further information on the
dynamics near $\cE_1$. To this end we prove the following result.
\begin{Thm}~\label{MP6}
  Let $f:\PP{2}\self$ be a holomorphic map of degree $d\ge2$ and
  let $S=\omega+dd^cu$ be a positive closed current on $\PP{2}$
  such that:
\begin{itemize}
  \item
    $S$ does not charge any component of $\cE_1$;
  \item
   the Lelong number $\nu(p,S)=0$ at any point $p\in\cE_2\cap\cE_1$.
  \end{itemize}
  Then
  \begin{equation}\label{Me11}
    \sup_{p\in\cE_1}\nu(p,d^{-n}f^{n*}S)\to0\quad\text{as $n\to\infty$}.
  \end{equation}
\end{Thm}

We first prove Theorem~\ref{MP6} under the weaker assumptions 
\begin{itemize}
\item[(A)]
  $u\not\equiv-\infty$ on any irreducible component of 
  $\cE_1$;
\item[(B)]
  $u$ is bounded at each point $p\in\cE_2\cap\cE_1$,
\end{itemize}

\begin{proof}[Proof of Theorem~\ref{MP6} under assumptions~(A) and~(B)]
  Let $V$ be an irre\-du\-ci\-ble component of $\cE_1$, i.e.\ a line.
  After replacing $f$ by an iterate we may assume that $V$ is fixed
  by $f$ and hence $R \= f|_V$ induces a rational map of $V$ of degree $d$.
  Let $S$ be a current as in the statement of the lemma. By~(A) we may
  define the probability measure $m_S \= S|_V$. We have
  \begin{equation*}
    \nu(p,d^{-n}f^{n*}S)
    \le \nu( p, d^{-n}R^{n*}m_S)
    = d^{-n} e ( p, R^n)\nu (p , m_S).
  \end{equation*}
  By~(B) the measure $m_S$ does not charge totally invariant orbits of $R$. On
  the other hand, one can find $\l<d$ such that 
  for large $n \ge 0$ and any $p\in V\setminus\cE_2$ 
  $e(p,R^n)\le\l^n$. We conclude 
  \begin{equation*}
    \sup_{p\in\cE_1}\nu(p,d^{-n}f^{n*}S)\le(\l/d)^n\to 0
    \quad\text{as $n\to\infty$}.
  \end{equation*}
\end{proof}
\begin{proof}[Proof of Theorem~\ref{MP6} in the general case.]
  We will use of Proposition~\ref{kisel} on 
  the behavior of Kiselman numbers when
  one weight tends to infinity.
  Let $V\subset\cE_1$ and $R\=f|_V$ be as above
  We cover $V$ by a finite number of coordinates chart $U_i \ni (z_i,
  w_i)$ such that $V \cap U_i = \{ z_i = 0 \}$. In the open set 
  $f^{-1}U_j\cap U_i$, the map $f$ can be written in the form
  \begin{equation*}
    f(z_i,w_i)=\left(z_i^d,R(w_i)+O(|z_i|)\right)\,(1+O(|z_i,w_i|)).
  \end{equation*}
  For a point $p \in U_i$, we denote by $\nu ( p, S, ( \a _1, \a _2))$ the
  Kiselman number of $S$ at $p$  with weights $ (\a _1, \a _2) \in (\R ^*
  _+)^2$ associated to the coordinate systems $(z_i, w_i)$. Assume that
  we can prove the following result:
  \begin{Lemma} \label{local}
    For any point $p \in f^{-1}U_j\cap U_i $ and any $\a \le 1$,
    we have:
    \begin{equation*}
      \nu\left(p,d^{-1}f^*S,(\a,d+1)\right)
      \le 
      \nu\left(fp,S,(\a\,e(p,R)/d,d+1)\right).
    \end{equation*}
  \end{Lemma}
  As before fix a constant $\l < d$ such that for large $ n \ge 0$ and
  any point $p \in V \setminus \cE_2$, we have  $ e ( p, R^n) \le \l ^n$.
  By assumptions $S$ does not charge $V$ and $\nu (p, S ) = 0$ for any
  point $p \in  V \cap \cE_2$; hence 
  for any $n\ge0$ large enough we get 
  \begin{align*}
    \nu(p,d^{-n}f^{n*}S)
    &\le C\,\nu\left(p,d^{-n}f^{n*}S,(1,d+1)\right)\\
    &\le C\,\nu\left(f^np,S,(e(p,R^n)/d^n,d+1)\right)\\
    &\le C\,\sup_{p\in V}\nu\left(p,S,(e(p,R^n)/d^n,d+1)\right)\\
    &\le C\,\sup_{p\in V}\nu\left(p,S,((\l/d)^n,d+1)\right)
    \longrightarrow_{n\to\infty}0,
  \end{align*}
  where the first inequality follows from~\eqref{Ec1} 
  (with $C\=d+1$), the second from Lemma~\ref{local}, 
  and the last convergence from Proposition~\ref{kisel}.
  This concludes the proof.
\end{proof}

\begin{proof}[Proof of Lemma~\ref{local}]
  This is a local result so we may assume $p=fp=(0,0)$  and
  \begin{equation*}
    f(z,w)=\left(
      z^d(1+O(|z|)),w^k\left(1+O(|w|)\right)+O(|z|)
    \right)\,
    (1+O(|z,w|)),
  \end{equation*}
  with $k=e(p,R)$. We
  easily check that there exist constants $C, C' > 0$ such that for any
  $\a\le 1$
  \begin{equation*}
    f\left(\Delta(r^{1/\a})\times\Delta(r^{1/(d+1)})\right)
    \supset
    \Delta(Cr^{d/\a})\times\Delta(C'r^{k/(d+1)}).
  \end{equation*}
  We remark that this property is easy to verify but nevertheless central
  to the proof.
  Write $S=dd^c u$ for some local psh potential $u$.
  We infer
  \begin{align*}
    \nu(p,d^{-1}f^{*}S,(\a,d+1))
    &=\lim_{r\to 0}\frac{\a(d+1)}{(\log r)}
    \sup_{\Delta(r^{1/(d+1)})\times\Delta(r^{1/\a})}u\circ f\\
    &\le\lim_{r\to 0}\frac{\a(d+1)}{(\log r)}
    \sup_{\Delta(Cr^{d/\a})\times\Delta(C'r^{k/(d+1)})}u\\
    &=k\,\nu\left(fp,S,(\a/d,(d+1)/k)\right)\\
    &=\nu\left(fp,S,(\a k/d,d+1)\right),
  \end{align*}
  which concludes the proof.
\end{proof}

%
%
%
%
\section{Proof of the main results}\label{proofmain}

This section is devoted to the proof of Theorem~A and its two
corollaries~B and~C.
\begin{proof}[Proof of Theorem~A]
  We argue by contradiction. Suppose that $S$ is a positive closed current
  on $\PP{2}$ for which the assumptions, but not the conclusions,
  of Theorem~A hold. As in Section~\ref{background} we 
  write $S=\omega+dd^cu$ with $u\le0$ qpsh, and conclude that
  there exists a ball $B$, a positive number $\alpha$ and a sequence
  $n_j\to\infty$ such that 
  \begin{equation}\label{Me5}
    f^{n_j}B\subset\{u<-\alpha d^{n_j}\}.
  \end{equation}
  We will get a contradiction from~\eqref{Me5}  by estimating
  the volumes of the two sides.
  
  Fix small neighborhoods $\Omega_1$, $\Omega_2$ of the exceptional sets
  $\cE_1$ and $\cE_2$, respectively. By the superattracting nature of
  $\cE_1$ and $\cE_2$ we may assume that $f\Omega_i\subset\subset\Omega_i$
  for $i=1,2$. In order to reach a contradiction, it is sufficient
  to consider three different cases.
  
  $\bullet$ Let us first assume that 
  $f^nB$ avoids $\Omega_1\cup\Omega_2$ for all $n\ge 0$. 
  Then Theorem~\ref{vol-t:1} applies and shows that
  \begin{equation*}
    \vol f^nB\ge(C_1\vol B)^{C_2 \lambda^n}
  \end{equation*}
  for some $\lambda<d$.
  On the other hand, the Kiselman-Skoda estimate (Theorem~\ref{vol-t:2})
  shows that 
  \begin{equation*}
    \vol\{u\le-\alpha d^n\}\le C\exp(-\beta d^n)
  \end{equation*}
  for some $\beta>0$ and for all $n\ge 0$. This yields a contradiction.
  
  $\bullet$ The second case is when $f^nB\subset\Omega_1\setminus\Omega_2$ 
  for all $n\ge 0$. We then use the results from 
  Sections~\ref{exc1},~\ref{volume1} and~\ref{lelong}
  on the dynamics near the first exceptional set $\cE_1$.
  First, by Proposition~\ref{MP4} there
  exists a constant $C>0$ such that
  \begin{equation}\label{Me12}
    \vol f^nB\ge(C\,\vol B)^{d^n}
  \end{equation}
  for sufficiently large $n$. Second, by Proposition~\ref{MP6}, for
  arbitrarily large $A >0$, one can
  find an integer $m\ge 0$ so that $\sup_{\cE_1} \nu ( p, d^{-m}
  f^{m*}S) <  1/A$. Hence
  by the Kiselman-Skoda estimate (Theorem~\ref{vol-t:2})
  one has
  \begin{equation}\label{Me13}
    \vol\{p\in\Omega_1\setminus\Omega_2\ \vert\
    d^{-m}u\circ f^m\le-t\}
    \le\exp(-A t)
  \end{equation}
  for large enough $t$.
  For $n_j \gg m $,~\eqref{Me5},~\eqref{Me12} and~\eqref{Me13} then imply
  \begin{align*}
    (C\,\vol B)^{d^{n_j -m}}
    &\le\vol f^{n_j-m}B\\
    &\le\vol\{ d^{-m} u\circ f^{m}<-\alpha d^{n_j - m} \}\\
    &\le\exp(-A \alpha {d^{n_j - m}}).
  \end{align*}
  We get a contradiction by choosing $A$ so that 
  $\exp(-\alpha A )<C\,\vol B$ and
  letting $n_j \to \infty$.
  
  $\bullet$ The third and last case is when 
  $f^nB\subset\Omega_2$ for all $n\ge0$.
  But by our assumption $u$ is bounded at $\cE_2$ and 
  so~\eqref{Me5} clearly cannot hold. 
  This completes the proof of Theorem~A.
\end{proof}  

\begin{proof}[Proof of Corollary~B]
  If $S=k^{-1}[C]$ is the current of integration
  on a curve $C$ of degree $k\ge1$, then $S$ satisfies the
  assumptions of Theorem~A unless
  \begin{itemize}
  \item
    $C$ contains an irreducible component of $\cE_1$; or 
  \item
    $C\cap\cE_2\ne\emptyset$.
  \end{itemize}
  This concludes the proof as the set of curves $C$ satisfying either
  of these conditions is a algebraic proper subset of $\PP{N}$.
\end{proof}

\begin{proof}[Proof of Corollary~C]
  Let ${\mathcal H}\subset\mathrm{Hol}_d$ be the set of holomorphic maps 
  $f$ of degree $d$ for which $\cE_f\ne\emptyset$.
  By Theorem~\ref{str-t:2} and Proposition~\ref{MP5}, 
  $\cE_f$ consists of at most
  three totally invariant lines and a  totally invariant set whose
  cardinality is bounded by some integer $N(d)$.
  It is to check from this 
  that ${\mathcal H}$ defines an algebraic set in $\mathrm{Hol}_d$.
  To conclude the proof we only have to exhibit one holomorphic map 
  $f\in\mathrm{Hol}_d$ with $\cE_f=\emptyset$. We follow a construction 
  of Ueda.
  
  Take a Latt\`es map in $\PP{1}$ of degree $d$ for instance 
  $R(z)\=(z-2/z)^d$. Consider the holomorphic map 
  $g(z,w)\=(R(z),R(w)):\PP{1}\times\PP{1}\self$ 
  It has topological degree $d^2$.
  The quotient $\PP{1}\times\PP{1}$ by the symmetry $(z,w)\to(w,z)$
  is isomorphic to $\PP{2}$ and $g$ induces a holomorphic map $f$ on the
  quotient. The topological  degree of $f$ is $d^2$ hence 
  $f\in\mathrm{Hol}_d$. As $R$ does not contain critical periodic points, 
  the same is true for $g$ and for $f$ too. Hence $\cE_f=\emptyset$ and we
  are done.
\end{proof}                                

%
%
%
%
\section{The proof of Theorem~A' and totally invariant currents}\label{misc}

In this section we will work under the assumption that 
\emph{every point in $\cE_2 \setminus \cE_1$ is a homogeneous point} 
i.e.\ $f$ preserves the pencil of
lines through that point. It is possible that this assumption is valid
for any holomorphic map of $\PP{2}$. Our goal is to prove Theorem~A'
and to exhibit totally invariant currents associated with the sets
$\cE_1$ and $\cE_2$.

%
%
\subsection{Local dynamics near $\cE_2$}\label{localform}

Near the points of $\cE_2$, the dynamics has a simple form and
this will allow us to prove good volume estimates.
\begin{Lemma}{~}\label{ML4}
  Assume $p$ is a homogeneous point. Then
  $f$ is locally conjugate at $p$ to a map of the form
  \begin{equation}\label{Me16}
    (z,w)\mapsto(P(z,w),Q(z,w)),
  \end{equation}
  where $P,Q$ are homogeneous polynomial of degree $d$.
\end{Lemma}
\begin{Lemma}{~}\label{ML5}
  Assume $p\in\cE_2\cap\cE_1$. Then 
  $f$ is locally conjugate at $p$ to a map of the form
  \begin{equation}\label{Me14}
    (z,w)\mapsto(z^d+wh(z,w),w^d),
  \end{equation}
  where $h$ is holomorphic.
\end{Lemma}

\begin{proof}[Proof of Lemma~\ref{ML4}.]
  Assume that $p=[0:0:1]$. In homogeneous
  coordinates, $f$ can be written 
  $f[z:w:t]=[P(z,w):Q(z,w):R(z,w,t)]$ for
  homogeneous polynomials $P$, $Q$, $R$ 
  of degree $d$ with $R(0,0,1)=1$.  
  Hence, locally, $f(z,w)=(P(z,w)(1+\eta),Q(z,w)(1+\eta))$
  for some germ $\eta$ with $\eta(0,0)=1$.
  As $f$ is contracting, one can define the map 
  \begin{equation*}
    \phi\=\prod_{j=0}^{\infty}(1+\eta\circ f^j)^{1/d^{j+1}}
  \end{equation*}
  and one checks the map
  $(z,w)\to\left(z\phi(z,w),w\phi(z,w)\right)$
  conjugates $f$ to $(P,Q)$.
\end{proof}
\begin{proof}[Proof of Lemma~\ref{ML5}.]
  Again assume $p=[0:0:1]$. We may assume that the set $\cE_1$ is
  given by $zw=0$ or by $w=0$. In the first of these cases, $p$
  is a homogeneous point and $f$ is locally conjugate to
  $(z^d,w^d)$ by Lemma~\ref{ML4}. In the second case, we have
  \begin{equation*}
    f[z:w:t]=[z^d+wQ(z,w,t):w^d:R(z,w,t)]
  \end{equation*}
  in homogeneous coordinates, where $R(0,0,1)=1$. Hence, locally,
  $f(z,w)=((z^d+wQ(z,w))(1+\eta),w^d(1 +\eta))$
  for some germ $\eta$ with $\eta(0,0)=1$.
  As in the proof of Lemma~\ref{ML4} we define
  \begin{equation*}
    \phi\=\prod_{j=0}^{\infty}(1+\eta\circ f^j)^{1/d^{j+1}}
  \end{equation*}
  and conclude that the map
  $(z,w)\to\left(z\phi(z,w),w\phi(z,w)\right)$
  conjugates $f$ to the desired form.
\end{proof}

\begin{Cor}\label{MC1}
  There exists $\alpha>0$ such that $c_n(p)\ge\alpha d^n$ for any
  $p\in\cE_2$.
\end{Cor}
\begin{proof}
  This follows immediately from the normal forms in 
  Lemma~\ref{ML4} and Lemma~\ref{ML5}.
\end{proof}

\begin{Prop}\label{MP8}
  Let $p\in\cE_2$ and let $\Omega$ be a small neighborhood of $p$.
  Then for any Borel set 
  $E\subset\Omega$ of positive volume $\vol (E) >0$, there exists $\gamma(E)>0$ such that
  \begin{equation}\label{Me17}
    \vol f^nE\ge \gamma(E)^{d^n}\quad\text{for all $n\ge0$}.
  \end{equation}
\end{Prop}
\begin{proof}
  We first consider the case $p\in\cE_2\cap\cE_1$ and write $f$ in the
  skew product form~\eqref{Me14}, which we may rewrite as
  \begin{equation}\label{Me15}
    (z,w)\mapsto(f_w(z),g(w))=(\psi\prod_{i=1}^d(z-\alpha_i(w)),w^d),
  \end{equation}
  where $\alpha_i$ are multi-valued with $\alpha_i(0)=0$ and
  $\psi(0,0)=1$. 

  Fix $\e_0>0$ small.
  It follows from~\eqref{Me15} that there exists a constant $c>0$
  such that for any $w\in\D(0,\e_0)$ and any Borel set 
  $E''\subset\D(0,\e_0)$ we have $\area f_wE''\ge c(\area E'')^d$.
  Further, for any Borel set $E'\subset\D(0,\e_0)$ we have
  $\area gE'\ge c(\area E')^d$. Iterating these estimates yields
  $\area g^nE'\ge (c'\area E')^{d^n}$ and
  $\area f^n_wE''\ge (c'\area E'')^{d^n}$, where
  $f_w^n=f_{g^{n-1}(w)}\circ\dots\circ f_w$.

  Now pick a Borel set $E\subset\Omega$ with $\vol E>0$. 
  After iterating forward
  we may assume that $E\subset\D^2(0,\e_0)$.
  For $w\in\D(0,\e_0)$ we write 
  $E_w''=\{z\in\D(0,\e_0)\vert\ (z,w)\in E\}$.
  There exists $\delta>0$ and a set 
  $E'\in\D(0,\e_0)$ with $\area E'>\delta$ such that
  $\area E_w''>\delta$ for $w\in E'$. 
  But then the previous estimates imply that
  $\area g^n E'>(c'\delta)^{d^n}$ and 
  $\area f_w^nE_w>(c'\delta)^{d^n}$ for $w\in E'$,
  so by Fubini's Theorem we get 
  $\vol f^nE\ge\gamma^{d^n}$ as desired.
  
  The remaining case, when $p\in\cE_2$ is a homogeneous point,
  is similar. We use the skew product structure~\eqref{Me16}.
  The only new observation that we need is that if $g:\PP{1}\self$
  is a rational map of degree $d\ge2$, then there exists $c>0$
  such that $\area g^nE\ge(c\area E)^{d^n}$ for any Borel set $E$.
\end{proof}

After these preliminaries we now prove Theorem~A'.
\begin{proof}[Proof of Theorem~A']
  The implication (1)$\Rightarrow$(2) is relatively easy. If $S$
  puts mass on a totally invariant curve $V\subset\cE_1$, say 
  $S\ge c[V]$, then for all $n\ge 0$ we have
  $d^{-n}f^{n*}S\ge d^{-n}f^{n*}c[V]=c[V]$.
  Since $T$ has bounded potential we cannot have convergence towards $T$.
  Similarly, if $p\in\cE_2$ with $\nu(p,S)>0$, then
  one immediately checks that $\nu(p,d^{-n}f^{n*}S)\ge d^{-n}c_n\nu(S,p)$.
  Hence, by Corollary~\ref{MC1},
  $\nu(p,d^{-n}f^{n*}S)\ge\e$ for some $\e>0$, which also prevents the
  sequence to converge towards $T$.
  
  Conversely, suppose that the current 
  $S$ satisfies~(2) of Theorem A'. 
  To prove that $d^{-n}f^{n*}S\to T$ we follow 
  the proof of Theorem~A up to the third case, i.e.\ 
  when $f^nB\subset\Omega_2$ 
  for all $n\ge0$. We pick a constant $\e\ll 1$ small enough.
  As $\nu(p,S)=0$ for all $p\in\cE_2$, by
  Theorem~\ref{vol-t:2} one can find a constant $C_{\e}>0$
  such that 
  \begin{equation}\label{Mis1}
    \vol\left\{u\le-t\right\}\le C_{\e}\exp(-t/\e),
  \end{equation}
  for all $t \ge 0$.
  Combining~\eqref{Mis1} and the hypothesis 
  $f^{n_j}B\subset\{u<-\alpha d^{n_j}\}$ 
  with Proposition~\ref{MP8}, we get
  \begin{equation*}
    \gamma(B)^{d^{n_j}}\le C_{\e}\exp(-\a d^{n_j}/\e) , 
  \end{equation*}
  hence $\gamma(B)\le\exp(-\a/\e)$ by letting $n_j\to+\infty$. But $\gamma(B)>0$
  is fixed and $\e$ is arbitrarily small, so this yields a contradiction.
\end{proof}

%
%
\subsection{Totally invariant currents}\label{exccurr}

Let us discuss the existence of totally invariant currents.
Consider the cone $\cS$ of positive closed
currents on $\PP{2}$ of unit mass such that $d^{-1}f^*T=T$.
Let $\cS^e$ be the set of extremal points in $\cS$. It is 
known~\cite{FoSi2} that the Green current $T$ is in $\cS^e$
(this follows e.g.\ from Theorem~A).

The following result is an immediate consequence of Theorem~A.
\begin{Cor}
  If $\cE_1=\cE_2=\emptyset$, then $\cS^e=\cS=\{T\}$.
\end{Cor}

Conversely we want to show that if either $\cE_1$ or $\cE_2$
is nonempty, then $\cS^e$ contains currents other than $T$.
Recall that the Green current $T$ has zero 
Lelong number at every point and, in particular,
does not put mass on any curve in $\PP{2}$. 
\begin{Prop}
  If the first exceptional set $\cE_1$ is nonempty, then there
  exists a current $S\in\cS^e$ supported on $\cE_1$.
\end{Prop}
\begin{proof}
  Since $\cE_1$ is totally invariant, the current 
  $(\deg\cE_1)^{-1}[\cE_1]$ is in $\cS$. This current 
  need not be extremal,
  but can be decomposed into currents in $\cS^e$ supported on $\cE_1$.
\end{proof}

In the sequel, $\|S\| \= \int S \wedge\omega$ denotes the projective mass of
the positive closed current $S$.

\begin{Prop}
  If $p\in\cE_2 \setminus \cE_1$ is a homogeneous point, 
  then there exists $S\in\cS^e$ with positive
  Lelong number at $p$ and with continuous potential outside $p$.
  More precisely, we have 
$ \nu (p, S ) = \| S \| = 1$.
\end{Prop}
\begin{proof}
  Assume  $f$ preserves
  the pencil of lines through $p$. Then $f$ induces a rational map
  $g$ of $\PP{1}$ (the set of lines) of degree $d$. Let $\mu$ be the measure
  of maximal entropy for $g$. This satisfies $g^*\mu=d\,\mu$ and 
  if we define
  \begin{equation*}
    S=\int_{a\in\PP{1}}[L_a]d\mu(a),
  \end{equation*}
  where $[L_a]$ denotes the current of integration on the line through
  $p$ corresponding to $a\in\PP{1}$, then $f^*S=d\,S$, so
  $S\in\cS$. Assume $S_1 \le S$ with $f^* S_1 = d S_1$. From the local
  structure of $S$, we infer the existence of a positive measure
  $\mu_1$ such that $  S_1=\int_{a\in\PP{1}}[L_a]d\mu_1(a)$. The
  equation $f^* S_1 = d S_1$ is equivalent to $f^* \mu_1 = d
  \mu_1$. As $\mu_1$ has no atoms, this
  forces $\mu_1 = c \mu$ for some constant $c>0$. Hence $ S_1 =
  c S$, and $S \in\cS^e$.
  Finally, a direct computation yields $\nu(p,S)=1$.
\end{proof}

\begin{Prop}
  If $p\in\cE_2 \cap \cE_1$,
  then there exists $S\in\cS^e$ with positive
  Lelong number at $p$ and with continuous potential outside $p$.

  More precisely, for any positive real numbers $0 < \a \leq 1$,
  there exists a current  $S_{\a}\in\cS^e$ with 
$$\nu (p,S_{\a},(1,\a / d^{-1}) ) = \| S_{\a}\| = 1.$$
\end{Prop}
\begin{proof}
Pick $\a \in ]0,1 ]$.
Introduce the homogeneous real analytic function on
$\CC{3}\setminus\{0\}$
$$
U(z,w,t) \= |z|^d+|w|^d+|w^{\a}t^{d-1}|.
$$
It vanishes exactly on the ray $\C \cdot(0,0,1)$ and by homogeneity the
positive closed current $d^{-1} dd^c \log U$ 
can be pushed down to $\PP{2}$ as a current
$\omega_0$, smooth outside $p$, with a pole at $p$ whose Lelong
number is $\nu(\omega_0,p) =1$.

Define $$
V(z,w,t) \= \frac{U \circ F}{U^d}
$$
where $F = (z^d+wQ,w^d,R)$ is a lift a $f$ to $\CC{3}$. Then $V$ induces by projection a
 function on $\PP{2}$ which is real analytic outside $p$. We claim there exists
constants $C_1, C_2 >0$ such that
$$
C_1 \le V \le C_2
$$
for any point in $\PP{2}$.

Indeed, as $p$ is totally invariant, $U$ and $U \circ F$ both vanishes
exactly along $\C \cdot (0,0,1)$, hence the inequality has to be
checked only in a neighborhood $V$ of $p$ in the chart $\{ t=1\}$. 
To do so, you may decompose $V$ in two sets $\{ | z^d| \leq A |w^{\a}|\}$ and 
$\{ | z^d| \geq A |w^{\a}|\}$ for a well chosen $A$.
In each of these sets, the estimates follow from a direct computation  
we leave to the reader. By normalizing $U$, we can
assume $C_2=1$.

We now follow the standard construction of the Green current. We have   
$f^* \omega_0 = d \omega_0 + dd^c \log V$, hence 
$$
d^{-k}f^{k*} \omega_0 = \omega_0 + dd^c\left( \sum_{j=0}^{k} d^{-j-1}\log
V \circ f^j\right)
$$
for all $k \ge0$.
The sequence of function $ \sum_{j=0}^{k} d^{-j-1}\log
V \circ f^j$ is decreasing converging uniformly on compact sets. Hence
the limit $G_p$ is a $L^1$ function, continuous  outside $p$, and
bounded everywhere.
The positive closed current $S \= \omega_0 + dd^c G_p$ belongs to $\cS$, has a
continuous potential outside $p$, and a singularity at $p$ with Lelong number
$\nu (p,S)=1/d$.
More precisely, in the coordinates $z,w$, the Kiselman number of $S$ with
weight $(1,\a/ d^{-1})$ is given by
$$
\nu (p,S,(1,\a/ d^{-1}) )=1.
$$
To conclude we show $S$ is extremal in the cone $\cS$. Assume $S =
S_1 + S_2$ with $f^* S_i = d S_i$ for $i = 1,2$. Kiselman numbers
behave additively hence
$$
\nu (p,S_1,(1,\a/d^{-1}) ) +  \nu (p,S_2,(1,\a/d^{-1}) ) = \nu
(p,S,(1,\a/d^{-1}) )=1.
$$
On the other hand, the following inequalities are standard (see \cite{Kis2},\cite{De1})
$$
\nu (p,S_i,(1,\a/d^{-1}) ) \le \nu(p,S_i) \le \| S_i\|.
$$

Whence
\begin{equation*}
1 =
\nu (p,S_1,(1,\a/d^{-1}) ) + \nu (p,S_2,(1,\a/d^{-1}) )
\leq
  \| S_1\| \cdot \| S_2\| =  \| S\| =1 
\end{equation*}
and we infer
\begin{equation}\label{EC1}
\nu (p,S_i,(1,\a/d^{-1}) ) =  \| S_i\|
\end{equation}
for $i=1,2$. Pick a global potential $S_i = dd^c u_i$ defined in
$\CC{2}$. By definition of Kiselman numbers we have in a neighborhood of
the origin
$$
u_i - \nu (p,S_i,(1,\a/d^{-1}) ) \times d^{-1} \log ( |z^d| + |w^{\a}| ) 
\le C
$$
for some constant $C>0$.
Together with equation \eqref{EC1}, we deduce that
$v_i \= u_i -  \| S_i\| d^{-1} \log ( |z^d| + |w^{\a}|  + |w^d| )$ is
globally bounded from above. As $u _1 + u_2 = G_p + d^{-1} \log (
|z^d| + |w^{\a}|  + |w^d| )$, and $G_p$ is bounded, we conclude that
$v_1,v_2$ are  also bounded
everywhere.
 
Hence
\begin{eqnarray*}
S_i
&=&
\lim_{n \to \infty} d^{-n} f^{n*} S_i 
\\
&=&
\lim_{n \to \infty}  dd^c ( d^{-n}  v_i \circ f^n ) +   \| S_i\| 
d^{-n} f^{n*} \omega_0
=
 \| S_i\| \cdot S
\end{eqnarray*}
This shows that $S \in \cS^e$.
\end{proof}

\begin{Example}
  For the map $f[z:w:t]=[z^d:w^d:t^d]$, the set $\cS^e$ is quite
  large. Given $\alpha,\beta,\gamma>0$ with 
  $\alpha+\beta+\gamma=1$ define
  \begin{equation*}
    u_{\alpha,\beta,\gamma}(z,w,t)
    =\alpha\log|z|+\beta\log|w|+\gamma\log|t|.
  \end{equation*}
  Then $\cS^e$ contains all the currents $S=\omega+dd^cu$ with
  \begin{equation*}
    u(z,w,t)=\max_{i=1,2}u_{\alpha_i,\beta_i,\gamma_i}(z,w,t)-\log|(z,w,t)|
  \end{equation*}
  such that
  \begin{equation*}
    \min_{i=1,2}{\alpha_i}
    =\min_{i=1,2}{\beta_i}
    =\min_{i=1,2}{\gamma_i}=0
  \end{equation*}
  or 
  \begin{equation*}
    u(z,w,t)=\max_{i=1,2,3}u_{\alpha_i,\beta_i,\gamma_i}(z,w,t)-\log|(z,w,t)|
  \end{equation*}
  such that
  \begin{equation*}
    \min_{i=1,2,3}{\alpha_i}
    =\min_{i=1,2,3}{\beta_i}
    =\min_{i=1,2,3}{\gamma_i}=0.
  \end{equation*}
  Notice that the Green current $T$ is of the latter form.
\end{Example}

%
%
\subsection{Configurations of exceptional sets.}\label{conf}

We conclude the paper by listing the different possible configurations of the
exceptional sets $\cE_1$ and $\cE_2$ and the corresponding mappings
$f$ in case $\cE_2 \setminus \cE_1$  contains only homogeneous points.
The case of totally invariant curves was treated in \cite{FoSi1}
(see Proposition~\ref{MP1}).  We summarize the results in
Table~\ref{table1} below.
\begin{itemize}
\item
  $P,Q,R$ denotes homogeneous polynomials in three variables $z,w,t$ 
  except if we state it otherwise;
\item
  $p_z$, $p_w$ and $p_t$ denote the points $[1:0:0]$, $[0:1:0]$ 
  and $[0:0:1]$, respectively;
\item
  $Z$, $W$ and $T$ denote the lines $(z=0)$, $(w=0)$ and $(t=0)$, 
  respectively;
\item
  in the cases where $\#\cE_1\ge2$ or $\#\cE_2\ge2$, 
  we mention only the maps preserving all irreducible components 
  of $\cE_1$ and each point in $\cE_2$. 
  To be complete, one has to add maps which permute these sets.
\end{itemize}
The proof is essentially elementary. There are essentially only two 
points to check: any intersection point between two irreducible
components of $\cE_1$ is in $\cE_2$; and 
if $\cE_2$ contains two homogeneous points 
$p,q$, then $\cE_1$ contains the line $H$ passing through $p$ and $q$.
The first of these statements is easy; for the second
note that $f^{-1}H$ is a union of lines passing through $p$ as $p$ 
is homogeneous, and also a union of lines passing through $q$. 
Hence $f^{-1}H=H$ is a totally invariant line. 
If we blow up $\PP{2}$ at the two points $p, q$, we can lift $f$ 
to a holomorphic map for which the strict transform of $H$ is 
totally invariant. We can hence contract it to a point, and the induced map 
becomes a holomorphic map on $\PP{1}\times\PP{1}$.
If $\{p,q\}=\{p_z,p_w\}$, this shows $f$ can be written under the form 
$f = [P(z,t):Q(w,t):t^d]$.
The other cases can be treated in a similar way.

In particular we have (assuming $\cE_2 \setminus \cE_1$
contains only homogeneous points):
\begin{Prop}
  There are at most $3$ distinct points in $\cE_2$.
\end{Prop}
\begin{table}
  \begin{tabular}{cc|ccc}
    $\#\cE_1$ & $\#\cE_2$ & $\cE_1$ & $\cE_2$ & $f$\\\hline
    1 & 0 & $T$     &               & $[P:Q:t^d]$\\
    0 & 1 &         & $p_w$         & $[P(z,t):Q:R(z,t)]$\\
    1 & 1 & $T$     & $p_z$         & $[P:w^d+tQ:t^d]$\\
      &   & $T$     & $p_t$         & $[P(z,w):Q(z,w):t^d]$\\    
    1 & 2 & $T$     & $p_z,p_w$     & $[P(z,t):w^d+tQ:t^d]$\\
    2 & 1 & $W,T$   & $p_z$         & $[P:w^d:t^d]$\\
    2 & 2 & $W,T$   & $p_z,p_w$     & $[z^d+tP:w^d:t^d]$\\
    2 & 3 & $W,T$   & $p_z,p_w,p_t$ & $[z^d+wtP:w^d:t^d]$\\
    3 & 3 & $Z,W,T$ & $p_z,p_w,p_t$ & $[z^d:w^d:t^d]$
  \end{tabular}

  \medskip
  \caption{Configuration of exceptional sets.}
  \label{table1}
\end{table}

%
%
%
%

%
%
%
%

\end{document}